\documentclass[12pt]{article}

\usepackage{graphics}
\usepackage{amsmath,amsfonts,latexsym,amssymb}
\usepackage{comment}
\usepackage{xcolor}
\usepackage{fullpage}
\usepackage{epsfig}  \input{psfig}

\def\Z{\mathbb{Z}}

\def\G{\Gamma}
\def\R{\mathcal{R}}
\def\C{\mathcal{C}}
\def\E{\mathcal{E}}
\def\S{\mathcal{S}}

\newtheorem{theorem}{Theorem}[section]
\newtheorem{lemma}[theorem]{Lemma}
\newtheorem{corollary}[theorem]{Corollary}
\newtheorem{conjecture}[theorem]{Conjecture}
\newtheorem{proposition}[theorem]{Proposition}
\newtheorem{definition}[theorem]{Definition}
\newtheorem{example}[theorem]{Example}

\newtheorem{remark}[theorem]{Remark}

\begin{document}
\title{A survey of Heffter arrays}

\author{ Anita Pasotti\thanks{Partially supported by INdAM - GNSAGA}
\\
DICATAM - Sez. Matematica\\ Università degli Studi di Brescia\\ Via Branze 43\\ I-25123 Brescia\ \  Italy\\
{\tt anita.pasotti@unibs.it} \\
\and
Jeffrey H. Dinitz \\
Dept of Math. and Stat.\\ University of Vermont\\ Burlington, VT 05405, USA\\
{\tt jeff.dinitz@uvm.edu}
}  

\maketitle              

\begin{center} \em
Dedicated to Doug Stinson on the occasion of his 66th birthday.
\end{center}

\begin{abstract}
Heffter arrays were introduced by Archdeacon in 2015 as an interesting link between combinatorial designs and topological graph theory.  Since the initial paper on this topic, there has been a good deal of interest in Heffter arrays as well as in related topics such as the sequencing of subsets of a group, biembeddings of cycle systems on a surface, and orthogonal cycle systems.
This survey presents an overview of the current state of the art of this topic. We begin with an introduction to Heffter arrays for the reader who is unfamiliar with the subject,  then we give a unified and comprehensive presentation of the major results, showing some proof methods also.  This survey also includes sections on the connections of Heffter arrays to several other combinatorial objects, such as problems on partial sums and sequenceability, biembedding graphs on surfaces, difference families and orthogonal graph decompositions. These sections are followed by  a section discussing the variants and generalizations of Heffter arrays which have been proposed.  The survey itself is complemented by a list of unsolved problems as well as an updated and complete bibliography.
\end{abstract}

\section{Introduction and Background}
The notion of a Heffter array was  introduced by Archdeacon \cite{A} in 2015.  The original motivation was to aid in an embedding problem in topological graph theory, but they were also introduced as a very interesting combinatorial design object.   The name is derived from
their connection to Heffter's first difference problem, a problem related
 to the construction of Steiner triple systems \cite{CR}.
Specifically, in 1896 Heffter \cite{H} introduced his famous \emph{First Difference Problem}
which asks  if, given $v\equiv 1 \pmod 6$, the set $\left\{1,2,\ldots,\frac{v-1}{2}\right\}$ can be
partitioned into $\frac{v-1}{6}$ triples $\{x, y, z\}$ such that either $x+y = z$ or $x+y+z=v$.
The problem was eventually solved in 1939 by Peltesohn \cite{P}.  We will show later in this survey that this problem is essentially solved (with an even stronger property) by a Heffter array with 3 rows and $\frac{v-1}{6}$ columns.

We begin with some definitions.
Let $\Z_v$ be the cyclic group of odd order $v$ whose elements are denoted by $0$ and $\pm i$ where
$i = 1, \ldots, \frac{v-1}{2}$. A \emph{half-set} $L$ of $\Z_v$ is a subset of $\Z_v\setminus \{0\}$ of size $\frac{v-1}{2}$
 that contains
exactly one of each pair $\{x,-x\}$. A \emph{Heffter system} $D(v, k)$, on a half-set $L$ of $\Z_v$, is a partition of $L$ into
parts of size $k$ such that the elements in each part sum to $0$ modulo $v$, see \cite{A,MR}.
Clearly a necessary condition for the existence of such a system is $v \equiv 1 \pmod{2k}$.

\begin{example}\label{ex:D315}
A Heffter system $D(31,5)$:
$$\{\{6,7,-10,-4,1\}, \{-9,5,2,-11,13\}, \{3,-12,8,15,-14\}\}.$$
\end{example}
It is easy to see that starting from a Heffter system $D(v,3)$ one immediately gets a solution to the Heffter's First
Difference Problem, as shown in the following example.

\begin{example}\label{ex:D313}
A Heffter system $D(31,3)$:
$$\{\{6,3,-9\}, \{7,5,-12\}, \{2,8,-10\},\{-4,-11,15\},\{1,13,-14\}\}.$$
\end{example}
Then $\{\{6,3,9\}, \{7,5,12\}, \{2,8,10\},\{4,11,15\},\{1,13,14\}\}$ is a solution to the Heffter's First
Difference Problem with $v=31$.  Using this Heffter system we get the following base blocks $\{\{0,6,9\}, \{0 ,7, 12\}, \{0, 2,10\},\{0,4,15\},\{0,1,14\}\}$ which can then be developed in $\Z_{31}$ to construct a Steiner triple system of order 31.  A detailed discussion of the construction of cycle systems from Heffter systems can be found in Section \ref{sec:DF}.

Two Heffter systems  $D_m = D(2mn + 1,m)$ and $D_n = D(2mn + 1, n)$ on the same
half-set $L$ of $\Z_{2mn+1}$ are \emph{orthogonal} if each part of $D_m$  intersects each
part of $D_n$  in exactly one
element, see \cite{A}. By a direct check one can see that the Heffter systems of Examples \ref{ex:D315} and \ref{ex:D313} are orthogonal.

 In \cite{A} Archdeacon defined a \emph{Heffter array} $H(m, n)$ as an $m \times n$ array whose rows form a $D(2mn + 1, n)$
and whose columns form a $D(2mn + 1,m)$;  these are called the row and column Heffter
systems, respectively. Note that a Heffter array $H(m, n)$ is equivalent to a pair of orthogonal Heffter
systems. An $H(m,n)$ gives two orthogonal Heffter systems by definition; its associated row system and its associated column system.
Conversely, two orthogonal Heffter systems $D(2mn + 1,m)$ and $D(2mn + 1, n)$ give rise to an $H(m,n)$
whose cell of indices $i,j$ contains the only element which the $i$-th part of the first Heffter system shares with the $j$-th part of the second one.

\begin{example}\label{ex:35}
Starting from the orthogonal Heffter systems of Examples \ref{ex:D315} and \ref{ex:D313}
we obtain the following Heffter array $H(3,5)$ whose elements belong to $\Z_{31}$:
$$
\begin{array}{|r|r|r|r|r|} \hline
6 & 7 & -10 & -4 & 1 \\ \hline
-9 & 5 & 2 & -11 & 13 \\ \hline
3 & -12 & 8 & 15 & -14 \\ \hline
\end{array}
$$
\end{example}

In the same paper Archdeacon proposed a variation of Heffter arrays  which allows for some empty cells. Two Heffter systems
$D_h=D(2mh+1, h)$ and $D_k=D(2nk+1, k)$ on the same half-set of size $mh = nk$ are \emph{sub-orthogonal} if
each part of $D_h$ intersects each part of $D_k$ in {\em at most} one element (see \cite{AHL} for information on sub-orthogonal  factorizations). We note here that the term {\em sub-orthogonal} has not been widely adopted, hence for  the remainder of this survey we will
 use  {\em orthogonal} to denote either orthogonal or sub-orthogonal. As
before, we can form an $m \times n$ array, denoted by $H(m, n; h,k)$, whose cell of indices $i,j$ either
contains the common element in the $i$-th part
of $D_h$ and the $j$-th part of $D_k$, if any, or is empty otherwise.
The $H(m,n)$ defined above is nothing but a $H(m,n;n,m)$.

In other words, when empty cells are allowed, a Heffter array can be defined as follows (this is now the standard definition):
\begin{definition}\label{def:Heffter}
 A \emph{Heffter array} $H(m,n;h,k)$ is an $m \times n$ matrix with entries from $\mathbb{Z}_{2nk+1}$ such that:
\begin{itemize}
  \item[{\rm (a)}] each row contains $h$ filled cells and each column contains $k$ filled cells;
\item[{\rm (b)}]  for every $x \in \mathbb{Z}_{2nk+1} \setminus \{0\}$, either $x$ or $-x$ appears in the array;
\item[{\rm (c)}] the elements in every row and column sum to $0$ in $\mathbb{Z}_{2nk+1}$.
\end{itemize}
\end{definition}

\begin{example}\label{ex61284}
An $H(6, 12; 8, 4)$ over $\Z_{97}$:
$$
\begin{array}{|r|r|r|r|r|r|r|r|r|r|r|r|}
\hline
&& -1 & 2 &   5 & -6 &&& -25 & 26 & 29 & -30 \\ \hline
&& 3 & -4 & -7 &   8 &&& 27 & -28 & -31 & 32 \\ \hline
9 & -10 &&& -13 & 14 & 33 & -34 &&& -37 & 38 \\ \hline
-11 & 12&&& 15 & -16 & -35 & 36 &&& 39 & -40 \\ \hline
-17 & 18 & 21 & -22 &&& -41 & 42 & 45 & -46 && \\ \hline
 19 & -20 & -23 & 24 &&& 43 &-44 & -47 & 48 && \\ \hline
\end{array}
$$
\end{example}
Note that by the pigeon-hole principle zero cannot appear in the array, hence
 Condition {\rm (b)} of Definition \ref{def:Heffter} is equivalent to requiring that
the elements appearing in the array form a half-set of  $\Z_{2nk+1}$.
It is easy to see that trivial necessary conditions for the existence of an $H(m, n; h,k)$ are $mh = nk$, $3\leq k \leq m$, and
$3\leq h \leq n$.

If a Heffter array $H(m, n; h,k)$ is square, that is $m = n$, then necessarily $h=k$. In this case
it is a $H(n, n; k,k)$ and is denoted by $H(n; k)$.  And just to clarify, a  $H(n, k)$  (note the comma in place of the semicolon) is a Heffter array with no empty cells and with $n$ rows and $k$ columns,
or a $H(n, k; k,n)$.

\begin{example}\label{ex:54}
An $H(5; 4)$ over $\Z_{41}$:
$$\begin{array}{|r|r|r|r|r|} \hline
 & 17 & -8 & -14&5  \\ \hline
1& & 18& -9 &-10\\ \hline
-6 & 2 &  &19 &-15\\ \hline
-11&-12&3&&20 \\ \hline
16&-7&-13&4& \\ \hline
\end{array}
$$
\end{example}

\begin{definition}\label{def:Heffter.integer}
A Heffter array is an \emph{integer Heffter array} if condition {\rm (c)} of Definition \ref{def:Heffter}
 is strengthened so that the elements in every row and every column, seen as integers in $\pm\{1,2,\ldots,nk\}$, sum to zero in $\Z$.
\end{definition}
Note that the arrays of Examples \ref{ex:35}, \ref{ex61284} and \ref{ex:54} are integer Heffter arrays.
In the  following example we present  a non-integer Heffter array.

\begin{example}
A non-integer $H(7;3)$ over $\Z_{43}$:
$$\begin{array}{|r|r|r|r|r|r|r|}
\hline
15&	-13&	-2&&&&	\\\hline							
-11&	14&		&-3&&&\\\hline							
-4&		&-8&	12&&&\\\hline							
  &-1&	10&	-9&&&	\\\hline		
&&&&								5&	21&	17\\\hline
&&&&								18&	6&	19\\\hline
&&&&								20&	16&	7\\\hline
\end{array}$$
\end{example}

A class of Heffter arrays of fundamental importance and which are particularly useful in recursive constructions of  are so called \emph{shiftable} arrays, see \cite{ADDY}.
An $n \times n$ array $A$ (possibly with empty cells) whose elements are integers is shiftable if each row and each column of $A$
contains the same number of positive and negative numbers.
Given a shiftable array $A$ and a nonnegative integer $x$, let $A \pm x$ denote the array  where $x$ is added to all the positive entries in $A$
and $-x$ is added to all the negative entries.
Define the {\em support} of $A$ as the set containing the absolute value of the elements appearing in $A$.
Note that if $A$ is shiftable with support $S$ and $x$ is a nonnegative integer, then $A\pm x$ has the same row and column sums as $A$ and has support $S + x$.
In the case of a shiftable integer Heffter array $H(m,n;h,k)$, the array
$H(m,n;h,k)\pm x$ has row and column sums equal to zero and its support is $\{1+x,2+x, \ldots ,nk+x\} $.
The arrays of Examples \ref{ex61284} and \ref{ex:54} are shiftable. The existence of shiftable square Heffter arrays was established
by constructive proofs in \cite{ADDY}.

\begin{theorem}{\rm\cite{ADDY}}\label{shift}
There exists a shiftable integer $H(n;k)$ if and only if $k$ is even and $nk\equiv 0 \pmod 4$.
\end{theorem}

Many of the known square Heffter arrays have a diagonal structure which has been shown to be useful for recursive constructions
and for the applications to biembeddings. In addition,  this diagonal structure has been extremely useful in computer searches for Heffter arrays.  Given a square array of order $n$, for $i=1,\ldots,n$, the $i$-th diagonal is  defined by
$D_i=\{(i,1),(i+1,2),\ldots, (i-1,n)\}$. All the arithmetic on row and column indices is performed modulo $n$, where the set of reduced
residues is $\{1,2,\ldots,n\}$. The diagonals $D_i,D_{i+1},\ldots, D_{i+k-1}$ are $k$ {\em consecutive diagonals}.
Given $n\geq k\geq 1$, a partially filled array $A$ of order $n$ is $k$-\emph{diagonal} if its nonempty cells are exactly those of $k$ diagonals.
Moreover, if these diagonals are consecutive, $A$ is said to be \emph{cyclically} $k$-\emph{diagonal}. The array of Example \ref{ex:54} is cyclically $4$-diagonal.

Shiftable cyclically $4$-diagonal $H(n;4)$ are constructed in \cite{ADDY} for all $n \ge 4$.  In view of their particular structure,
these arrays can be used to add four filled cells per row and column to an existing Heffter array $H$ if $H$
contains four consecutive diagonals of empty cells. This is described in the following lemma and example.

\begin{lemma}{\rm\cite{ADDY}}\label{lemma:shift}
If there exists an integer Heffter array $H(n;k)$ which  has $s$ disjoint sets of four consecutive empty diagonals,  then there exists an
integer
$H(n;k+4s)$.  Furthermore, if the $H(n;k)$ is shiftable, then the $H(n;k+4s)$ is shiftable too.
\end{lemma}

\begin{example}\label{ex:H87}
The following are a cyclically $3$-diagonal $H(8;3)$, say $H$, and a shiftable  $H(8;4)$, say $A$:
\begin{center}
\begin{footnotesize}
$H=\begin{array}{|r|r|r|r|r|r|r|r|}\hline
8 & 16 &   &  &  &  & & -24 \\
\hline  -17 & -6 & 23 &  &  &  &  &   \\
\hline   & -10 & -5 & 15 &   &  &  &  \\
\hline   &   & -18 & 7 & 11 &  &  &  \\
\hline   &  &  & -22 & 3 & 19 &   &  \\
\hline   &  &  &   & -14 & 2 & 12 &  \\
\hline    &  &  &  &  & -21 & 1 & 20 \\
\hline  9 &  &  &  &  &   & -13 & 4  \\\hline
\end{array}$
\end{footnotesize}

\vspace{1em}
\begin{footnotesize}
A= $\begin{array}{|r|r|r|r|r|r|r|r|}\hline
& &   & 1 & -3 & -5 & 7 &  \\
\hline   &  &  & -4 & 2 & 8 & -6 &   \\
\hline  15 &  &  &  &   & 9 & -11 & -13 \\
\hline  -14 &   &  &  &  & -12 & 10 & 16 \\
\hline  -19 & -21 & 23 &  &  &  &   & 17 \\
\hline  18 & 24 & -22 &   &  &  &  & -20 \\
\hline    & 25 & -27 & -29 & 31 &  &  &  \\
\hline   & -28 & 26 & 32 & -30 &   &  &  \\\hline
\end{array}$
\end{footnotesize}
\end{center}

Using $H$ and $A\pm 24$ one can construct the following $H(8;7)$ over $\Z_{113}$:
\begin{center}
\begin{footnotesize}
$\begin{array}{|r|r|r|r|r|r|r|r|}\hline
8 & 16 &   & 25 & -27 & -29 & 31 & -24 \\
\hline  -17 & -6 & 23 & -28 & 26 & 32 & -30 &   \\
\hline  39 & -10 & -5 & 15 &   & 33 & -35 & -37 \\
\hline  -38 &   & -18 & 7 & 11 & -36 & 34 & 40 \\
\hline  -43 & -45 & 47 & -22 & 3 & 19 &   & 41 \\
\hline  42 & 48 & -46 &   & -14 & 2 & 12 & -44 \\
\hline    & 49 & -51 & -53 & 55 & -21 & 1 & 20 \\
\hline  9 & -52 & 50 & 56 & -54 &   & -13 & 4  \\\hline
\end{array}$
\end{footnotesize}
\end{center}
\end{example}

Concerning the search for Heffter arrays constructed from diagonals, a useful notation was introduced  in \cite{DW} and has been used in several subsequent papers including  \cite{weak,RelH,CPEJC,CPPBiembeddings}.
All the constructions in \cite{DW}  were based on filling in the cells of a fixed collection of diagonals. To aid in the constructions the authors defined the following procedure for filling a sequence of cells on a diagonal. It is termed {\em diag} and it has six parameters, as follows.

In an $n \times n$ array $A$ the procedure
$diag(r,c,s,\Delta_1,\Delta_2, \ell)$ installs the entries of an array $A$ as follows:
\begin{center}
$A[r+i\Delta_1,c+i\Delta_1]=s+i\Delta_2\mbox{ \ for \ }i=0,1,\dots,\ell-1.$\end{center}
Here all arithmetic on the row and column indices is performed modulo $n$,
where the set of reduced residues is $\{1,2, \ldots, n\}$.
The following summarizes the parameters used in the $diag$ procedure:

\begin{tabular}{cl}
$\bullet$ & $r$ denotes the starting row,\\
$\bullet$ &$c$ denotes the starting column,\\
$\bullet$ & $s$ denotes the starting symbol,\\
$\bullet$ &  $\Delta_1$ denotes how much the row and column indices change at each step,\\
$\bullet$ & $\Delta_2$ denotes how much the symbol changes at each step, and\\
$\bullet$ & $\ell$ is the length of the chain.\\
\end{tabular}

\noindent
The following example shows the use of the above definition.

\begin{example}\label{ex1} {\rm\cite{DW}}
The array below is a shiftable integer Heffter array $H(11;4)$ where the filled cells are contained in two sets of consecutive diagonals
and it is
constructed via the following procedures:

\noindent

\begin{center}
$\begin{array}{ll}
  diag(4,1,1,1,2,11);\quad\quad & diag(5,1,-2,1,-2,11); \\
  diag(4,7,-23,1,-2,11);\quad \quad  & diag(5,7,24,1,2,11).
\end{array}$
\end{center}

\renewcommand{\arraycolsep}{1.5pt}
\begin{footnotesize}
\[
\begin{array}{|r|r|r|r|r|r|r|r|r|r|r|}\hline
 & &38&-39& & & &-16&17& & \\ \hline
 &&&40&-41&&&&-18&19&  \\ \hline
 &&&&42&-43&&&&-20&21 \\ \hline
1&& &&&44&-23&&&&-22 \\ \hline
-2&3&&&&&24&-25&&&  \\ \hline
 &-4&5&&&&&26&-27&&  \\ \hline
 &&-6&7&&&&&28&-29&  \\ \hline
 &&&-8&9&&&&&30&-31 \\ \hline
-33&&&&-10&11&&&&&32 \\ \hline
34&-35&&&&-12&13&&&&  \\ \hline
 &36&-37& & & &-14&15& && \\ \hline
\end{array}
\] \end{footnotesize}
\renewcommand{\arraycolsep}{6pt}
\end{example}

\medskip

Concerning the existence problem,
the first classes of Heffter arrays which were systematically studied were (1) square  arrays and  (2) Heffter arrays with no empty cells (termed {\em tight} Heffter arrays). Square integer Heffter arrays were considered in
 \cite{ADDY}. In that paper the authors proved the following almost complete result.
 In their proofs Lemma \ref{lemma:shift} plays a fundamental role.

\begin{theorem}{\rm\cite{ADDY}}
If an integer $H(n;k)$ exists, then necessarily $3 \leq k \leq n$ and $nk \equiv 0,3\pmod 4$.
Furthermore, these conditions are sufficient except possibly when $n\equiv 0\ \mbox{or} \ 3\pmod 4$ and  $k\equiv 1\pmod 4$.
\end{theorem}

It should be noted that \cite{ADDY} also contains partial results
when $n\equiv 0\ \mbox{or} \ 3\pmod 4$ and $k\equiv 1$ (mod
$4$).
These missing cases were entirely solved in \cite{DW} using the $diag$ method described above.
Thus we have  a complete constructive solution for the integer square case.
This result is summarized
 in the following theorem.

\begin{theorem}\label{thm:intergerHeffter}{\rm \cite{ADDY,DW}}
There exists an integer $H(n; k)$ if and only if $3 \leq k \leq  n$ and $nk \equiv
0, 3\pmod4$.
\end{theorem}

For all congruences classes of $n$ and $k$ for which an integer Heffter array necessarily cannot exist a non-integer Heffter array $H(n;k)$ was proven to exist in \cite{CDDY}.
Hence we can state the following.
\begin{theorem}\label{thm:Heffter}{\rm \cite{ADDY,CDDY,DW}}
There exists an $H(n; k)$ if and only if $3 \leq k \leq  n$.
\end{theorem}

Table \ref{table} details the papers where existence is proven for each congruence class of $n$ and $k$ modulo 4.
\begin{table}[h]
 \begin{center}
\begin{tabular} {|c||c|c|c|c|} \hline
$n\backslash k$&          0&                   1&          2&3            \\ \hline\hline
0              &\cite{ADDY}&\cite{ADDY,DW}&\cite{ADDY}&\cite{ADDY} \\ \hline
1              &\cite{ADDY}&                 \cite{CDDY}&       \cite{CDDY} &\cite{ADDY} \\ \hline
2              &\cite{ADDY}&                 \cite{CDDY} &\cite{ADDY}&     \cite{CDDY} \\ \hline
3              &\cite{ADDY}&\cite{ADDY,DW}&        \cite{CDDY}&        \cite{CDDY} \\ \hline
\end{tabular}
\end{center}
\label{table}
\caption{Articles which prove the existence of square Heffter arrays $H(n;k)$.}
\end{table}

We now turn to the case of
rectangular Heffter arrays without empty cells $H(m,n)$, or so-called tight Heffter arrays.
A complete constructive solution is presented in \cite{ABD}.  Historically speaking, it was this paper that introduced the notion of shiftable Heffter arrays. The main theorem of that paper is given next.

\begin{theorem}{\rm \cite{ABD}}
There exists an $H(m,n)$ if and only if $m,n\geq 3$.
Moreover, these arrays are integer if and only if $mn\equiv 0,3 \pmod 4$.
\end{theorem}

It has not been proven that integer rectangular arrays with empty cells exist for all possible orders. We note here that the necesary conditions for the existence of  an integer $H ( m, n ; h, k )$ are that
 $3 \leq h \leq n$, $3 \leq  k \leq m$, $mh = nk$ and $nk \equiv 0,3 \pmod 4$.
 Recent partial results for integer rectangular arrays with empty cells can be found in \cite{MP3}.  We summarize these results in the following two theorems.

	\begin{theorem}{\rm \cite{MP3}}\label{Fiore1}
Let $m, n, h, k$ be  such that $3 \leq h \leq n$, $3 \leq  k \leq m$ and $mh = nk$.
Set $d = \gcd ( h , k )$. There exists an integer
$H ( m, n ; h, k )$ in each of the following cases:
\begin{itemize}
\item[$(\rm{1})$] $d \equiv 0 \pmod 4$;
\item [$(\rm{2})$] $d \equiv 1 \pmod 4$ with $d \geq 5$ and $nk \equiv 3 \pmod 4$;
\item[$(\rm{3})$] $d \equiv 2 \pmod 4$ and $nk \equiv 0 \pmod 4$;
\item[$(\rm{4})$] $d \equiv 3 \pmod 4$ and $nk \equiv 0, 3 \pmod 4$.
\end{itemize}
\end{theorem}

\begin{theorem}{\rm \cite{MP3}}\label{Fiore2}
	Let $m, n, h, k$ be  such that $3 \leq h \leq n$, $3 \leq  k \leq m$ and $mh = nk$.
 If $h\equiv 0 \pmod 4$ and $k \neq 5$ is odd, then
there exists an integer  $H ( m, n ; h, k )$.
\end{theorem}

This survey is organized as follows.  In Section \ref{sec:simple} we discuss results on simple and globally simple Heffter arrays.  Then in Section \ref{sec:biemb} we give the application of Heffter arrays to biembedding complete graphs on topological surfaces.  In Section \ref{sec:sums} we consider the problem of sequencing subsets of a group.  Section \ref{sec:DF} gives results concerning the connection between Heffter arrays, difference families  and graph decompositions and in Section \ref{sec:variants} we survey results concerning variants and generalizations of Heffter arrays. Finally in Section \ref{sec:conclusions} we give some open problems related to Heffter arrays and the other objects covered in this survey.

\section{Simple Heffter Arrays}\label{sec:simple}

In \cite{A}, Archdeacon showed that Heffter arrays can be used to construct (simple) cycle decompositions of the complete graph
if they satisfy an additional condition, called \emph{simplicity}.
This concept is also related to several open problems on partial sums and sequenceability, these questions will be discussed in the Section \ref{sec:sums}.

Let $T$ be a finite subset of an additive group $G$. Given an ordering $\omega=(t_1,t_2,\ldots,t_k)$ of the elements in $T$, let $s_i=\sum_{j=1}^i t_j$
be the $i$-th partial sum of $\omega$.
The ordering $\omega$ is said to be \emph{simple} if its partial sums $s_1,\ldots,s_k$ are pairwise distinct.
We point out that, if $\sum_{t \in T}t=0$, then an ordering $\omega$ of the elements in $T$ is simple if and only if it does not contain any
proper subsequence that sums to $0$. We also note that if $\omega=(t_1,t_2,\ldots,t_k)$ is a simple ordering so is $\omega^{-1}=(t_k,t_{k-1},\ldots,t_1)$.

In the following, with a little abuse of notation, we will identify each row  and column of an $H(m,n;h,k)$ with the subset of $\Z_{2nk+1}$ of size $h$  and $k$, respectively,
whose elements are the entries of the nonempty cells of such a row  or column. For example we can view the first
row of the $H(8;7)$ of Example \ref{ex:H87} as the subset $R_1=\{8,16,25,-27,-29,31,-24\}$ of $\Z_{113}$.

\begin{definition}
An $H(m,n;h,k)$ is  \emph{simple} if each row and each column admits a simple ordering.  
\end{definition}

It is easy to see that since each row and each column of an $H(m,n;h,k)$  does not contain $0$ or $2$-subsets of the form
$\{x,-x\}$, and is such that $s_h=0$ and $s_k=0$, respectively, then every $H(m,n;h,k)$ with $h,k\leq 5$ is simple (indeed {\em any} subset of size $k\leq 10$ with these properties is simple \cite{AL}).

\begin{example}\label{ex:simple}
The $H(8;7)$ of {\rm Example \ref{ex:H87}} is simple. To verify this property we need to provide a simple ordering for each row and each column.
One can check that the $\omega_i$'s are simple orderings of the rows and the  $\nu_i$'s are simple orderings of the columns:
\begin{center}
\begin{footnotesize}
$\begin{array}{lcl}
\omega_1= ( 8, 25, 16, -27, -29, 31, -24); & \quad &\nu_1= (8,39,-17,-38,-43,42,9  );\\
\omega_2= (-17, -6, -28, 23, 26, 32, -30); & \quad &\nu_2= ( 16,-6,-45,-10,48,49,-52 );\\
\omega_3= (39, -10, -5, 33, 15, -35, -37); & \quad & \nu_3= ( 23,-5,47,-18,-46,-51,50 );\\
\omega_4= (-38, -18, 7, -36, 11, 34,  40); & \quad &\nu_4= (25,-28,15,7,-53,-22,56 );\\
\omega_5= (-43, -45, 47, -22, 3, 41,  19); & \quad &\nu_5= (-27,26,11,3,55,-14,-54 );\\
\omega_6= ( 42, 48, -46, -14, 2, -44, 12); & \quad &\nu_6= ( -21,32,33,-36,19,2,-29);\\
\omega_7= (20, -51, -53, 55, -21, 1, 49);& \quad  &\nu_7= (-13,-30,-35,34,12,1,31);\\
\omega_8= (-52, 9, 50, 56, -54, -13, 4); & \quad &\nu_8= (-37,-24,40,41,-44,20,4).
\end{array}$
\end{footnotesize}
\end{center}
\end{example}

By \emph{natural ordering} of a row (column, respectively) one means the ordering from left to right (from top to bottom, respectively).
Since the natural ordering of each row and each column of the $H(8;7)$ of Example \ref{ex:H87} contains a proper subsequence that sums to zero
(this immediately follows by its construction),
the natural ordering cannot be simple for any row and any column.
Clearly,  for larger $n$ and $k$ it is more difficult and tedious to provide
explicit simple orderings for rows and columns of an $H(m,n;h,k)$.
For this reason, in \cite{CMPPHeffter} the authors considered Heffter arrays which are simple with respect to the natural ordering
of rows and columns
and proposed the following definition.

\begin{definition}{\rm \cite{CMPPHeffter}}
A Heffter array $H(m, n;h, k)$ is \emph{globally simple} if the natural ordering of each row and each column is simple.
\end{definition}

A globally simple $H(m,n; h,k)$ is denoted by $SH(m,n; h,k)$ and a square globally simple $H(n;k)$ is
denoted by $SH(n;k)$.

By above arguments, we have that the $H(8;7)$ of Example \ref{ex:H87} is simple but not globally simple.  However  a globally simple $H(8;7)$ does indeed exist and is presented in the following example.

\begin{example}\label{ex:SH87}
A globally simple $H(8;7)$:
\begin{center}
\begin{footnotesize}
\begin{tabular}{|r|r|r|r|r|r|r|r|}
\hline $4$ & $35$ & $-45$ & $46$ & $ $ & $20$ & $-36$ & $-24$ \\
\hline $48$ & $-5$ & $23$ & $-47$ & $-18$ & $ $ & $37$ & $-38$ \\
\hline $-32$ & $-10$ & $-6$ & $31$ & $-41$ & $42$ & $ $ & $16$ \\
\hline $33$ & $-34$ & $44$ & $3$ & $11$ & $-43$ & $-14$ & $ $ \\
\hline $ $ & $15$ & $-28$ & $-22$ & $7$ & $27$ & $-53$ & $54$ \\
\hline $-13$ & $ $ & $29$ & $-30$ & $56$ & $1$ & $12$ & $-55$ \\
\hline $-49$ & $50$ & $ $ & $19$ & $-40$ & $-21$ & $2$ & $39$ \\
\hline $9$ & $-51$ & $-17$ & $ $ & $25$ & $-26$ & $52$ & $8$ \\\hline
\end{tabular}
\end{footnotesize}
\end{center}
\end{example}

The Heffter arrays constructed in {\rm \cite{ADDY,DW}}, in general, are not globally simple and
 easy modifications of the existing constructions seem not to produce globally simple arrays. One exception to this are the
$SH(n;6)$'s constructed in  {\rm \cite[Proposition 4.1]{CMPPHeffter}}
which were obtained by switching the first two columns of the matrices given in {\rm \cite[Theorem 2.1]{ADDY}}.
Globally simple Heffter arrays are more elegant than simple ones but also more difficult to construct,
in fact there are few papers on these arrays, see \cite{BCDY,CDY,CMPPHeffter,DM}.   The following two theorems summarize the results on square globally simple integer Heffter arrays. The first deals with the case when $k\leq 10$ and the second considers the existence of several of the congruence classes.

\begin{theorem}{\rm \cite{CMPPHeffter}}\label{thm:SH}
Let $3\leq k\leq 10$. There exists an integer globally simple $H(n;k)$ if and only if $n \geq k$ and $nk\equiv 0,3 \pmod4$.
\end{theorem}

\
\begin{theorem}{\rm \cite{BCDY}} \label{thm:globally_simple2}
Let $n\geq k\geq 3$. There exists an integer globally simple  $H(n;k)$ in each of the following cases:
\begin{itemize}
  \item[$(\rm{1})$ ] $k\equiv 0\pmod 4$;
  \item[$(\rm{2})$ ] $n\equiv 1\pmod4$ and $k\equiv3\pmod 4$;
  \item[$(\rm{3})$ ] $n\equiv 0\pmod4$, $k\equiv3\pmod 4$ and $n>>k$.
\end{itemize}
\end{theorem}

In the rectangular case we have the following theorem for Heffter arrays with either three or five rows.
\begin{theorem}{\rm \cite{DM}} \label{thm:SH1}
 There exists a globally simple $H(3,n)$ for every $n\geq 3$ and a globally simple $H(5,n)$ for all $3\leq n\leq 100$.
\end{theorem}

\section{Applications to biembeddings}\label{sec:biemb}

One of the main reasons that Archdeacon introduced the notion of Heffter arrays is their connection with biembedding cycle systems on a topological surface.
Several papers focus on this application, see, for instance \cite{CDY,CMPPHeffter,CPArxiv,DM}.
We emphasize that there has been previous work done in the area of biembedding cycle systems, especially  triple systems, see \cite{DGLM,DGLM2,GG,GK}.  It is interesting to
note that some of the earliest examples of Heffter arrays are derived from direct
constructions of biembeddings.  In particular, in \cite{ADDY} the authors used current graphs (developed by Youngs \cite{Y}) based on M\"{o}bius ladders with $n=4m+1$ rungs and on cylindrical ladders with $n=4m$ rungs to construct diagonal Heffter arrays $H(n;3)$. So not only do Heffter arrays yield biembeddings (as will be detailed in this section), but biembeddings have been used to make Heffter arrays.

We recall some basic definitions, see \cite{Moh} and also note that  \cite{GT} and \cite{MT} are excellent references for graphs on surfaces.
An \emph{embedding} of a graph (or multigraph) $\G$ in a surface $\Sigma$ is a continuous injective
mapping $\psi: \G \to \Sigma$, where $\G$ is viewed with the usual topology as $1$-dimensional simplicial complex.
The  connected components of $\Sigma \setminus \psi(\G)$ are called $\psi$-\emph{faces}.
If each $\psi$-face is homeomorphic to an open disc, then the embedding $\psi$ is said to be \emph{cellular}.
In this context, we say that two embeddings $\psi:\G\rightarrow\Sigma$ and $\psi':\G'\rightarrow\Sigma'$ are isomorphic if and
only if there is a graph isomorphism $\sigma:\G\rightarrow\G'$ such that $\sigma(F)$ is a $\psi'$-face if and only
if $F$ is a $\psi$-face.
A circuit decomposition of a graph $\G$ is a collection of edge-disjoint subgraphs of $\G$ such that every edge of $\G$ belongs to exactly one circuit.

\begin{definition}
A \emph{biembedding} of  two circuit decompositions $\mathcal{D}$ and $\mathcal{D}'$ of a simple graph $\G$  is a face $2$-colorable embedding
  of $\G$ in which one color class is comprised of the circuits in $\mathcal{D}$ and
  the other class contains exactly the circuits in $\mathcal{D}'$.
\end{definition}

	If the faces are cycles the biembedding is said to be \emph{simple}.
	When one is not interested in the circuits, we simply speak of a biembedding of the graph $\G$.
	Equivalently, a biembedding of a simple graph $\G$ is a face $2$-colorable embedding of $\G$.

\begin{example}\label{emb.fig}
A biembedding of the complete graph $K_7$ on the torus.  Note that the black cells are a Steiner triple system (a circuit decomposition of $K_7$ into triples) and the white cells are also a Steiner triple system.

\centerline{\hbox{
\psfig{figure=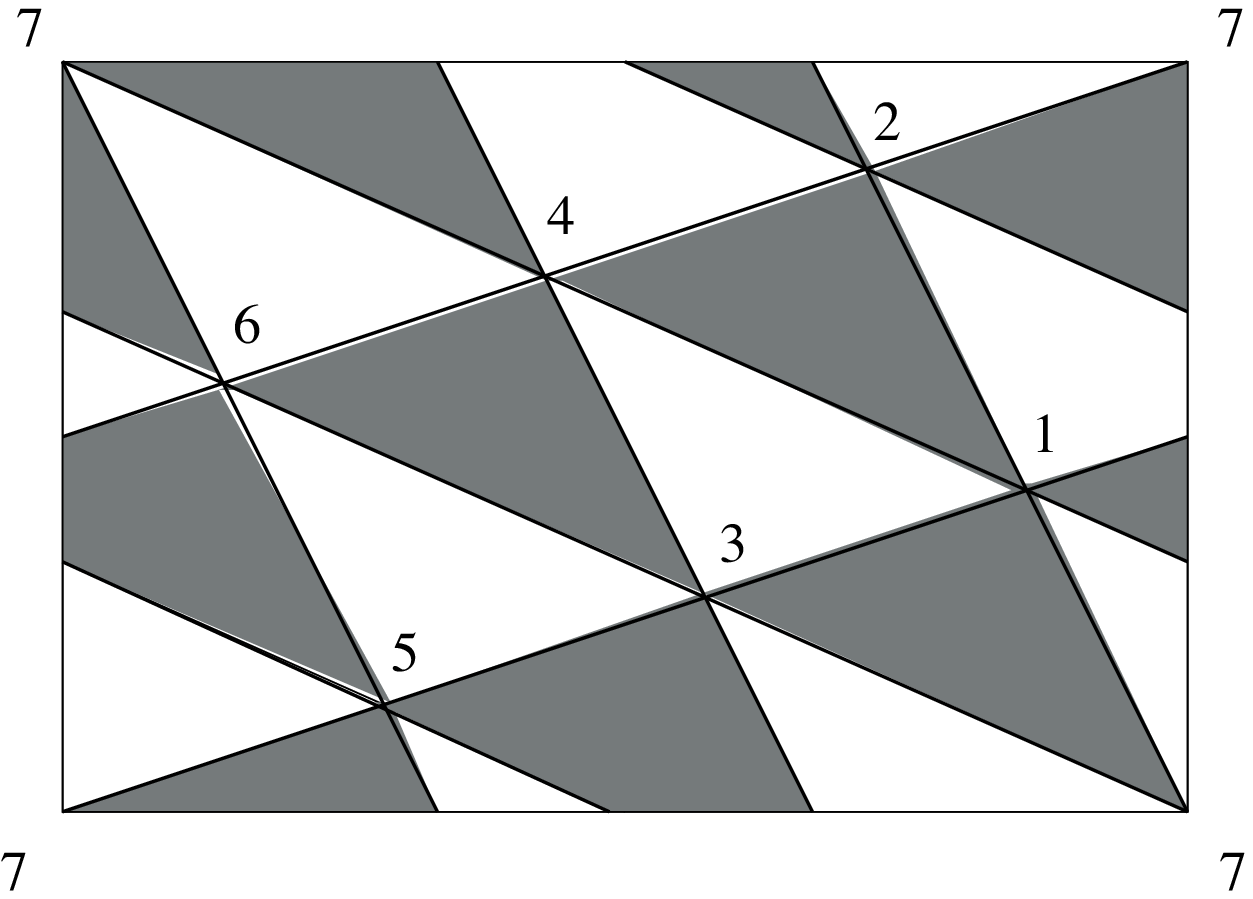,height=2in,width=3in}}}
\end{example}

In order to state our main theorem connecting Heffter arrays to biembeddings (Theorem \ref{Gustin}), we need the notion of compatible orderings of the rows and  columns of the Heffter array.

\begin{definition} Let $A$ be a partially  filled array with $F$ filled cells whose elements are pairwise distinct.
Let $\omega_r$ ($\omega_c$, respectively) be any ordering of the rows (columns, respectively) of $A$.
Then $\omega_r$ and $\omega_c$ are  \emph{compatible} orderings if $\omega_r\circ \omega_c$ is a cycle of order $F$.
\end{definition}

\begin{example}  Given the  $H(3,5)$ from Example \ref{ex:35}, use the natural ordering of the rows and columns to get row cycles
 $\omega_r= (6, 7, {-10}, {-4}, 1)(-9, 5, 2, {-11}, 13)(3, {-12}, 8, 15,{-14})$
 and the column cycles
 $\omega_c= (6,{-9},3)(7,5,{-12})({-10},2,8)({-4},{-11},15)(1,13,{-14}) $.
Then we see that
$\omega_r\circ \omega_c = ( 6,5,8, {-4}, 13, 3,7, 2,15,1,{-9},{-12},{-10},{-11},{-14}) $
is a single cycle (of order 15) and hence $\omega_r$ and $\omega_c$ are compatible orderings.
\end{example}

We will discuss the connection between compatible orderings and biembeddings below.  It is because of this connection that the existence of compatible orderings for a given array has been considered in several papers, see \cite{CDY,CDP,CMPPHeffter,DM}. The following theorem details the necessary conditions for there to exist
 orderings $\omega_r$ and $\omega_c$ of the row and column Heffter systems of a Heffter array that are compatible.

\begin{theorem}{\rm \cite{CDY,DM}}\label{thm:compatible}
If there exist compatible orderings $\omega_r$ and $\omega_c$ for a Heffter array
$H (m, n; h, k)$, then 
either:
\begin{itemize}
\item[{\rm (1)}] $m, n, h$ and $k$ are all odd;
\item[{\rm (2)}] $m$ is odd, $n$ is even and $h$ is even; or
\item[{\rm (3)}] $m$ is even, $n$ is odd and $k$ is even.
\end{itemize}
\end{theorem}

Hence, if  a square Heffter array $H(n;k)$ admits compatible orderings, then $nk$ is odd.
This condition is also sufficient when $k$ is small and the array is cyclically $k$-diagonal.
In fact the following holds:

\begin{proposition}{\rm \cite{CDP}}
Let $3\leq k< 200$ be an integer and
let $A$ be a cyclically $k$-diagonal Heffter
array of size $n\geq k$. Then there exist two compatible orderings $\omega_r$ and $\omega_c$ of
the rows and the columns of $A$ if and only if $n$ and $k$ are both odd.
\end{proposition}

In view of Theorem \ref{thm:compatible},  in order to obtain results on biembeddings, several authors have  focused on
square Heffter arrays with $n$ and $k$ odd.
The following results pertain to the existence of compatible orderings in the case of cyclically $k$-diagonal Heffter arrays $H(n;k)$.

\begin{proposition}{\rm \cite{CDP,CMPPHeffter}}
 Let $k\geq 3$ be an odd integer and let $A$ be a cyclically $k$-diagonal Heffter
 array of size $n\geq k$. If $\gcd(n, k-1) = 1$, then there exist two compatible orderings $\omega_r$ and $\omega_c$ of
the rows and the columns of $A$.
\end{proposition}

\begin{proposition}{\rm \cite{CDP}}
Let $k\geq3$  and let $A$ be a cyclically $k$-diagonal Heffter
array of  size $n\geq (k-2)(k-1)$. Then there exist two compatible orderings $\omega_r$ and $\omega_c$ of
the rows and the columns of $A$ if and only if $n$ and $k$ are both odd.
\end{proposition}

In the  case of an $H(m,n;h,k)$ with no empty cells (an $H(m,n)$) the following theorem was proven in \cite{DM}. In the proof, an explicit compatible ordering of the rows and columns is given.

\begin{proposition}{\rm \cite{DM}} \label{compat}
Let $A$ be a $m\times n$ Heffter array with no empty cells where at least one of $m$ and $n$ is odd. Then there exist compatible orderings, $\omega_r$ and $\omega_c$ on the rows and columns of $A$.
\end{proposition}

So from Theorem \ref{thm:compatible} and Proposition \ref{compat} we have necesary and sufficient condition for an $H(m,n)$ to have compatible orderings.  Specifically we have the following.

\begin{theorem}
Let $A$ be a $m\times n$ Heffter array with no empty cells (an $H(m,n))$.  Then there exist compatible orderings, $\omega_r$ and $\omega_c$ on the rows and columns of $A$ if and only if at least one of $m$ and $n$ is odd.
\end{theorem}

Compatible orderings are the link between Heffter arrays and biembeddings, while the simplicity of a Heffter array will guarantee the biembedding is simple also.  The following fundamental theorem from \cite{A} gives the connection between Heffter arrays with compatible orderings and biembeddings.  The interested reader is referred to that paper for the full details of this connection, but briefly, using the compatible orderings of the Heffter array one obtains  a current graph which in turn gives the rotation of the edges around each vertex.  The rows and columns give the cycles in the two cycle decompositions.

\begin{theorem}\label{Gustin}{\rm \cite{A}}
Given a Heffter array $H(m,n;h,k)$ with  compatible orderings $\omega_r$ on the rows and $\omega_c$ on the columns,
there exists a biembedding of the complete graph $K_{2nk+1}$
 on an orientable surface such that every edge is
on a face of size $h$ and a face of size $k$. Moreover, if $\omega_r$ and $\omega_c$ are both simple, then all
faces are simple cycles.
\end{theorem}

A $k$-cycle system of order $n$ is a set of cycles of length $k$, denoted by $C_k$, whose edges partition the edges of $K_n$,
one speaks also of a $C_k$-decomposition of  $K_n$.
Restating the previous theorem in terms of biembeddings of cycle systems we have the following.
\begin{theorem}\label{biembed2}{\rm \cite{A}}
Given a Heffter array $H(m,n;h,k)$ with simple compatible orderings $\omega_r$ on the rows and $\omega_c$ on the columns, there exists an orientable biembedding
of a $C_h$-decomposition and a $C_k$-decomposition both of $K_{2nk+1}$, or equivalently, there is an orientable biembedding
of an $h$-cycle system and a $k$-cycle system both in $K_{2nk+1}$.
\end{theorem}

We note again that a discussion of the construction of cycle systems from Heffter arrays will be given in Section \ref{sec:DF}.
As an example of the application of Proposition \ref{compat} and Theorems \ref{Gustin} and \ref{biembed2}  simple $H(3,n)$ were explicitly constructed in \cite{DM} for all $n\geq 3$ with an eye towards the biembedding problem.
The following is the main result in that paper.

\begin{theorem}{\rm \cite{DM}}
There exists a globally simple $H(3,n)$ for every $n\geq 3$.  Thus for every $n \geq 3$ there exists a biembedding of $K_{6n+1}$ on an orientable surface such that each edge is on a $3$-cycle and an $n$-cycle and  hence there is  a biembedding of a Steiner triple system and an $n$-cycle system.
\end{theorem}

The following problem was introduced in \cite{CDP} in view of its importance in looking for compatible orderings for Heffter arrays.
However, it can also be considered as a tour problem on a toroidal board which is interesting \emph{per se}.

Let $A$ be an $m\times n$ toroidal partially filled  array. By $r_i$ we denote the orientation of the $i$-th row of $A$,
precisely $r_i=1$ if it is from left to right and $r_i=-1$ if it is from right to left. Analogously, for the $j$-th
column, if its orientation $c_j$ is  from  top to bottom then $c_j=1$ otherwise $c_j=-1$. Assume that an orientation
$\R=(r_1,\dots,r_m)$
 and $\C=(c_1,\dots,c_n)$ is fixed. Given an initial filled cell $(i_1,j_1)$ consider the sequence
$ L_{\R,\C}(i_1,j_1)=((i_1,j_1),(i_2,j_2),\ldots,(i_\ell,j_\ell),$ $(i_{\ell+1},j_{\ell+1}),\ldots)$
where $j_{\ell+1}$ is the column index of the filled cell $(i_\ell,j_{\ell+1})$ of the row $R_{i_\ell}$ next
to
$(i_\ell,j_\ell)$ in the orientation $r_{i_\ell}$,
and where $i_{\ell+1}$ is the row index of the filled cell of the column $C_{j_{\ell+1}}$ next to
$(i_\ell,j_{\ell+1})$ in the orientation $c_{j_{\ell+1}}$.
The problem is the following:\\

\textbf{Crazy Knight's Tour Problem.}{\rm \cite{CDP}}
Given a toroidal partially filled array $A$,
do there exist $\R$ and $\C$ such that the list $L_{\R,\C}$ covers all the filled
cells of $A$?\\

The Crazy Knight's Tour Problem for a given array $A$ is denoted by $P(A)$.
Also, given a filled cell $(i,j)$, if $L_{\R,\C}(i,j)$ covers all the filled positions of $A$, then
 $(\R,\C)$ is said to be a solution of $P(A)$.

\begin{example}\label{ex0}
Let $A$ be the following array, where the bullets  ``$\bullet$'' denote the filled cells:
$$\begin{array}{|c|c|c|c|}\hline
     \bullet &  & & \bullet      \\\hline
      &   \bullet &   \bullet &      \\\hline
      &   \bullet  &  \bullet & \bullet    \\\hline
      \bullet &    &    \bullet  &    \\\hline
    \end{array}$$
		Choosing $\R:=(-1,1,1,-1)$ and $\C:=(1,-1,1,1)$ we can cover all the filled positions of $A$, as shown in the table below where
		in each filled cell we write $j$ if this position is the $j$-th element of the list $L_{\R,\C}(1,1)$.
The elements of $\R$ and $\C$ are represented by an arrow.		
$$\begin{array}{r|r|r|r|r|}
  &  \downarrow   & \uparrow & \downarrow & \downarrow   \\\hline
\leftarrow &     1 &      &    &  5  \\\hline
\rightarrow &      &   3 &   7 &      \\\hline
\rightarrow &      &   8  & 4  & 2   \\\hline
\leftarrow &      6 &    &  9  &    \\\hline
    \end{array}$$
\end{example}

\begin{remark}
 If $A$ is an $H(m, n; h, k)$ such that $P(A)$ admits a solution $(\R, \C)$, then $A$
admits two compatible orderings $\omega_r$ and $\omega_c$ that can be determined as follows.
For each row $R$ (column $C$, respectively) of $A$, let $\omega_R$ be
the natural ordering if $r_i = 1$ ($c_i = 1$, respectively) and its inverse otherwise.
Setting $\omega_r = \omega_{R_1} \circ \ldots \circ \omega_{R_m}$ and
$\omega_c = \omega_{C_1} \circ \ldots \circ \omega_{C_n}$, since $(\R, \C)$ is a solution of $P(A)$,
then  $\omega_r$ and
$\omega_c$ are compatible.
\end{remark}

The relationship between globally simple Heffter arrays, Crazy Knight's
Tour Problem and biembeddings is explained in the following result.  It  is a reformulation of Theorem \ref{Gustin}
in the case of square globally simple Heffter arrays in terms of $P(A)$. Clearly, a
similar theorem holds in the rectangular case.
\begin{theorem}
Let $A$ be a globally simple Heffter array $H(n;k).$  If there exists a solution of $P(A)$, then there exists a biembedding of
 $K_{2nk + 1}$ on an orientable surface whose face boundaries are $k$-cycles.
\end{theorem}

The existence of  integer globally simple  Heffter arrays $H(n;k)$ with compatible orderings for $n\equiv 1\pmod 4$ and $k\equiv 3 \pmod 4$ has been
proved (with infinite sporadic exceptions) in \cite{CDY}.

\begin{theorem}\label{HeffterBiembeddingsCavenagh}{\rm \cite{CDY}}
  Let $n \equiv 1\pmod 4$, $t > 0$ and $n > 4t + 3$. If there exists $\alpha$ such that
  $2t + 2\leq \alpha \leq  n-2-2t$, $\gcd(n, \alpha) = 1$, $\gcd(n, \alpha-2t-1) = 1$ and $\gcd(n, n-1-\alpha-2t) =
1$, then there exists a globally simple integer Heffter array $H (n; 4t + 3)$ with row and column orderings that are
both simple and compatible.
\end{theorem}

\begin{theorem}{\rm \cite{CDY}}
  Let $n \equiv 1\pmod 4$,  $n > k\equiv 3 \pmod 4$ and either: $(a)$ $n$ is prime; $(b)$
$n = k + 2$; or $(c)$  $n\geq 7(k + 1)/3$ and if $n \equiv 3 \pmod 6$ then $k \equiv 7\pmod {12}$. Then there
exists a globally simple integer Heffter array $H (n; k)$ with an ordering that is both simple
and compatible. Furthermore [by Theorem \ref{biembed2}], there exists a face $2$-colorable embedding $\mathcal{G}$
of $K_{2nk+1}$ on an orientable surface such that the faces of each color are cycles of length $k$.
Moreover, $\Z_{2nk+1}$ has a sharply vertex-transitive action on $\mathcal{G}$.
\end{theorem}

\begin{definition}Two Heffter arrays $H$ and $H'$ are said to be {\em equivalent} if one can
be obtained from the other by (i) rearranging rows or columns; (ii) replacing every entry $x$
in $H$ with $-x$; and/or (iii) taking the transpose.
\end{definition}

The following theorem gives a lower bound on the number of inequivalent Heffter arrays
$H(n; 4t + 3)$ that satisfy Theorem \ref{HeffterBiembeddingsCavenagh}.
Let $\mathcal{H}(n)$ represent
the number of derangements on $\{0,1,\ldots,n-1\}$.
\begin{theorem}{\rm \cite{CDY}}\label{BoundCavenagh}
 Let $n \equiv 1\pmod 4$, $t \geq 2$ and $n > 4t + 3$. If there exists $\alpha$ such that
$2t + 2\leq\alpha \leq  n-2-2t$, $\gcd(n, \alpha) = 1$, $\gcd(n, \alpha-2t-1) = 1$ and $\gcd(n, n-1-\alpha-2t) = 1$,
then there exist at least $(n -2)(\mathcal{H}(t-2))^2\simeq (n-2)[(t-2)!/e]^2$ inequivalent
globally simple integer Heffter arrays $H(n; 4t + 3)$, each with an ordering that is both simple
and compatible.
\end{theorem}

The other case, that is $n\equiv 3\pmod 4$ and $k\equiv 1 \pmod 4$, remains unsolved in general.
In the same paper the authors also investigated the number of biembeddings of $K_v$ obtained from certain classes of Heffter arrays.
 The following result is a corollary of Theorem \ref{BoundCavenagh}.

\begin{theorem}{\rm \cite{CDY}}
 Let $n \equiv 1\pmod 4$ and $k \equiv 3\pmod 4$, $k \geq 11$ and either $n$ is prime or
$n \geq (7k + 1)/3$. Further if $n \equiv 3\pmod 6$ assume $k \equiv 7\pmod {12}$. Then there exists at
least $(n-2)[((k-11)/4)!/e]^2$ distinct face $2$-colorable embeddings of the complete graph
$K_{2nk+1}$ onto an orientable surface where each face is a cycle of fixed length $k$, and the vertices
can be labeled with the elements of $\Z_{2nk+1}$ in such a way that this group $(\Z_{2nk+1}, +)$ has a
sharply vertex-transitive action on the embedding.
\end{theorem}

Starting from the results of \cite{CDY}, in \cite{CPArxiv} the authors consider the number of non-isomorphic biembeddings of $K_v$ obtained from certain classes of Heffter arrays.
Set the binary entropy function by $H(p):=-p\log_2{p}-(1-p)\log_2(1-p)$.

\begin{theorem}\label{GeneralBound}{\rm \cite{CPArxiv}}
Let $v=2nk+1$, $k=4t+3$ and let $n\equiv 1 \pmod{4}$ be either a prime or $n\geq (7k+1)/3$.
Moreover, if $n\equiv 3\pmod{6}$ assume $k\equiv 7\pmod{12}$. Then, the number of non-isomorphic
face $2$-colorable embeddings of the complete graph
$K_{2nk+1}$ onto an orientable surface where each face is a cycle of fixed length $k$
 is at least
$$\frac{(n-2)[\mathcal{H}(t-2)]^2}{2(2nk)^2}\approx \frac{\pi(t-2)^{2t-5}}{64e^{2t-2}n}\approx k^{k/2+o(k)}/v.$$
\end{theorem}

\begin{theorem}{\rm \cite{CPArxiv}}
  Let $k\geq 3$ be odd,  $n\geq 4k-3$ and  $nk\equiv 3 \pmod{4}$. Assume also that $\gcd(n,k-1)=1$. Then, set $v=2nk+1$, the number of non-isomorphic face $2$-colorable embeddings of the complete graph
$K_{2nk+1}$ onto an orientable surface where each face has fixed length $k$
 is at least:
$$ \frac{\binom{\lceil n/(k-1)\rceil}{ \lceil n/(4k-4)\rceil}}{(nk)^2}\gtrsim \frac{\sqrt{\frac{ 2(k-1)}{3n\pi}}2^{\frac{n}{k-1}\cdot H(1/4)+1}}{(nk)^2}\approx 2^{v\cdot \frac{H(1/4)}{2k(k-1)}+o(v,k)}.$$
Moreover, if:
\begin{itemize}
\item[$(\rm{1})$ ] $k=3$ the number of \emph{simple}  non-isomorphic face $2$-colorable embeddings of $K_{6n+1}$ is at least $\frac{2^{n/2}}{9n^2} \approx 2^{v/12+o(v)}$;\\
\item[$(\rm{2})$ ] $k=5,7,9$ the number of \emph{simple}  non-isomorphic face $2$-colorable embeddings of $K_{2nk+1}$ is at least $ \frac{\binom{\lfloor n/(k-1)\rfloor}{ \lfloor n/(4k-4)\rfloor}}{(nk)^2}\gtrsim \frac{\sqrt{\frac{ 2(k-1)}{3n\pi}}2^{\lfloor\frac{n}{k-1}\rfloor\cdot H(1/4)+1}}{(nk)^2}\approx 2^{v\cdot \frac{H(1/4)}{2k(k-1)}+o(v,k)}.$
\end{itemize}
\end{theorem}

One additional embedding result follows directly by Theorem \ref{thm:SH}.
\begin{theorem}{\rm \cite{CMPPHeffter}}
There exists a biembedding of the complete graph of order $2nk + 1$
and one of the cocktail graph of order $2nk +2$ on orientable surfaces such that every
face is a $k$-cycle, whenever $k \in \{3, 5, 7, 9\}$, $nk \equiv 3 \pmod 4$ and $n > k$.
\end{theorem}

\section{Sequencing Subsets of a Group}\label{sec:sums}

We know from Theorem \ref{Gustin} that if there is a simple Heffter array $H(m,n;h,k)$ (with compatible orderings), then there is a biembedding of $K_{2nk+1}$
 on an orientable surface such that every edge is
on a face of size $h$ and a face of size $k$ and all
faces are simple cycles.  It is thus a natural question to ask just when the rows and columns of a Heffter array have a simple ordering, or more generally,  when subsets of a group have a simple ordering.  Indeed, if this is always the case, then  Heffter arrays would always be simple and the condition for simplicity of the array would be unnecessary.  It has indeed been an active area of research to show this, but at this time it is still not known for all subsets (and all groups).
In this section we give some conjectures relating to the sequencing of subsets of groups as well as some of the results in this area.

We begin with the more classical notion of a sequenceable group.
Let $G$ be a group of order $n$ (written additively). Given an ordering $\omega=(0,g_1,g_2,\ldots,g_{n-1})$ of the elements in $G$, let $s_i=\sum_{j=1}^i g_j$
be the $i$-th partial sum of $\omega$. The group $G$ is {\em sequenceable} if there exists some ordering of the elements of $G$ where all of these partial sums
$s_0=0,s_1,\ldots,s_{n-1}$ are distinct
(equivalently, if $G$ admits a simple ordering), while it is {\em R-sequenceable} if $s_0=0,s_1,\ldots,s_{n-2}$ are distinct
and in addition $s_{n-1} = 0$.
The concept of sequenceable group was introduced in 1961 by Gordon \cite{G},
even if similar ideas for cyclic groups were already presented in 1892 by Lucas \cite{L}.
Alspach et al. \cite{AKP} proved that any finite abelian group is either sequenceable or R-sequenceable,
confirming the Friedlander-Gordon-Miller conjecture.

Here is a short summary of the results concerning the sequencing and R-sequencing of groups.  For an extensive survey on this topic the reader can refer to \cite{Osurvey}.

\begin{proposition}Results on sequencing of groups.

\begin{enumerate}
  \item[$(\rm{1})$ ]  An abelian group is sequenceable if and only if it has a unique element of order $2$.

\item[$(\rm{2})$ ]  No non-abelian group of order less than $10$ is sequenceable. But Keedwell conjectures that all non-abelian groups of order greater than $10$ are sequenceable.

\item[$(\rm{3})$ ]  Keedwell's conjecture has been proved true for the following classes of groups:
\begin{enumerate}
\item[$(\rm{a})$ ]  All non-abelian groups of order $n$, $10\le n\le 32$;
\item[$(\rm{b})$ ]  The dihedral groups $D_n$;
\item[$(\rm{c})$ ]  $A_5$ and $S_5$;
\item[$(\rm{d})$ ]  Solvable groups with a unique element of order $2$;
\item[$(\rm{e})$ ]  Some non-solvable groups with a unique element of order $2$;
\item[$(\rm{f})$ ]  Some groups of order $pq$, $p$ and $q$ odd primes.
\end{enumerate}
\end{enumerate}

\end{proposition}

\begin{proposition} Results on R-sequencing of groups.

\begin{enumerate}
\item[$(\rm{1})$ ]  If $G$ is an R-sequenceable group, then its Sylow $2$-subgroup is either trivial or non-cyclic.

\item[$(\rm{2})$ ]  The following groups are R-sequenceable:
\begin{enumerate}
\item[$(\rm{a})$ ]  $\Z_n$, $n$ odd;
\item[$(\rm{b})$ ]  Abelian groups of order $n$, $gcd(n,6)=1$;
\item[$(\rm{c})$ ]  $D_n$, $n$ even;
\item[$(\rm{d})$ ]  $Q_{2n}$ if $n+1$ is a prime of the form $4k+1$, for which $-2$ is a primitive root;
\item[$(\rm{e})$ ]  Non-abelian groups of order $pq$, $p<q$ odd primes, with $2$ a primitive root of $p$.
\end{enumerate}
\end{enumerate}
\end{proposition}

The problem of ordering the elements of a {\em subset} of a given group in such
a way that all the partial sums are pairwise distinct is more recent than the sequenceable group question.  Nonetheless it has also been extensively studied, most recently because of its connection to Heffter arrays, but not only for this reason. Indeed there are several open conjectures on partial sums of a finite group.

For example, several years ago  Alspach made
 the following conjecture, whose validity would shorten some cases of known proofs
about the existence of cycle decompositions.

\begin{conjecture}\label{Conj:als}
Let $T\subseteq \mathbb{Z}_v\setminus\{0\}$ such that
$\sum_{t\in T}t\neq0$. Then there exists an ordering of the elements of $T$ such that the partial sums are all
distinct and non-zero.
\end{conjecture}

The first published paper on this problem is by Bode and Harborth \cite{BH} in which the authors proved that Conjecture \ref{Conj:als}
is valid if $|T|=v-1$ or $|T|=v-2$ or if $v\leq 16$ (the latter was obtained by computer verification).
Recently, several other papers on this topic have been published, see \cite{AL,CDF,CDFORF,CP,HOS,O}.

Another conjecture, very close to the Alspach's conjecture, was originally proposed in 1971 by Graham \cite{Gr71}
for cyclic groups of prime order and then extended in 2015 by
 Archdeacon et al.  \cite{ADMS}  to every cyclic group as follows.
\begin{conjecture}\label{Conj:ADMS}{\rm \cite{ADMS}}
Let $T\subseteq \Z_v\setminus\{0\}$.
Then  $T$ admits a simple ordering.
\end{conjecture}

Clearly Conjecture  \ref{Conj:ADMS} implies Conjecture \ref{Conj:als}.  In \cite{ADMS} the authors proved that Conjecture \ref{Conj:als} implies Conjecture \ref{Conj:ADMS}.
Then they proved that Conjecture \ref{Conj:ADMS} is valid if $|T|\leq 6$ and (by computer verification) for $v\leq 25$.

In \cite{CMPPSums}, Costa et al. proposed the following conjecture on partial sums  which was also motivated by the study of Heffter arrays.

\begin{conjecture}\label{Conj:nostra}{\rm \cite{CMPPSums}}
Let $(G,+)$ be an abelian group. Let $T$ be a finite subset of $G\setminus\{0\}$
such that no $2$-subset $\{x,-x\}$ is contained in $T$ and
with the property that $\sum_{t\in T} t=0$.
Then  $T$ admits a simple ordering.
\end{conjecture}

If $G=\Z_v$, then every row and column of a Heffter array can be viewed as a subset $T$ considered in
Conjecture \ref{Conj:nostra}. Hence if this conjecture were true for cyclic groups, then every Heffter array would be simple.
In \cite{CMPPSums} the validity of Conjecture \ref{Conj:nostra} is proved theoretically for subsets $T$ of size less than $10$
and, with the aid of a  computer, for every abelian
group $G$ with $|G| \leq 27$.

Clearly if $G=\Z_v$, Conjecture \ref{Conj:nostra} immediately follows from Conjecture \ref{Conj:ADMS}.
Indeed in \cite{CMPPSums} the authors extended Conjecture \ref{Conj:ADMS} to every abelian group
and proved by computer that this extended conjecture  is valid for every abelian group of order at most $23$.
Also, they pointed out that in \cite{ADMS}, when proving the conjecture true for $|T|\leq 6$ the authors do not use
the hypothesis that $T$ is a subset of a cyclic group. In fact, these proofs for $|T|\leq 6$ work in general for an arbitrary abelian group.

Some results are known about Conjectures \ref{Conj:als} and \ref{Conj:ADMS} on an arbitrary group.
For instance the validity of Conjecture \ref{Conj:als} has been verified by computer
for \emph{all abelian groups} of order at most $21$  in \cite{CMPPSums}.  It has also been verified in \cite{AL}
for any subset $T$ of size at most $9$ of an \emph{arbitrary abelian group}.
On the other hand, Conjecture \ref{Conj:als} cannot be generalized to nonabelian groups as shown in \cite{An} by Anderson.
Also, if we consider Conjectures \ref{Conj:als} and \ref{Conj:ADMS} in the case of abelian groups then, again,
Conjecture \ref{Conj:als} implies Conjecture \ref{Conj:ADMS}. Hence it follows that
Conjecture \ref{Conj:ADMS} holds, in the abelian case, for any subset $T$ of size at most $9$ and  when $\sum_{t\in T}t=0$ for  $T$ of size
at most $10$, for details see \cite{ADMS}. This implies that also Conjecture \ref{Conj:nostra} holds for
 subsets of size at most $10$.
For Conjecture \ref{Conj:ADMS} it is also natural to investigate the nonabelian case.
In \cite{CMPPSums} this conjecture is proved theoretically for subsets $T$ of size at most $5$ of \emph{every arbitrary group}
and a computer verification is made for \emph{all groups} of order not exceeding $19$.

Finally, we recall that a subset of an (additive) abelian group is said to be \emph{sequenceable} if there is an ordering $(t_1,\ldots,t_k)$
of its elements such that the partial sums $(s_0,s_1,\ldots,s_k)$, where  $s_0 = 0$, are distinct, with the possible exception that we may
have $s_k=s_0=0$.
The following theorem summarizes what is known about sequencing subsets of the cyclic group.
We should note that most of the recent progress on this problem has been made by use of the {\em polynomial method}.  This
method relies on the non-vanishing corollary to the Combinatorial Nullstellensatz (see  \cite{Alon,CDF,CDFORF,HOS,O}).

\begin{theorem}
Let $T \subseteq \Z_v\setminus \{0\}$ with $|T|=k$.  Then $T$ is sequenceable in the following cases:
\begin{itemize}
\item[$(\rm{1})$ ] $k\leq 9$\  {\rm \cite{AL}};
\item[$(\rm{2})$ ] $k \leq 12$ when $v=pq$ where $p$ is prime and $q\leq 4$ \  {\rm \cite{CDFORF,HOS}};
\item[$(\rm{3})$ ] $k \leq 12$ when $v=mq$ where all the prime factors of $m$ are bigger than $\frac{k!}{2}$ and $q \leq 4$ \  {\rm \cite{CDFORF}};
\item[$(\rm{4})$ ] $k=v-3$ when $v$ is prime and $\sum_{t \in T}t=0$\  {\rm \cite{HOS}};
\item[$(\rm{5})$ ] $k=v-2$ when $\sum_{t \in T}t\neq 0$ \  {\rm \cite{BH}};
\item[$(\rm{6})$ ] $k=v-1$ \  {\rm \cite{AKP,G}};
\item[$(\rm{7})$ ] $v\leq 21$ and $v\leq 25$ when $\sum_{t \in T}t=0$ \  \rm{ \cite{ADMS,CMPPSums}}.
\end{itemize}
\end{theorem}

\section{Connections with difference families and graph decompositions}\label{sec:DF}

In this section we discuss the connection between Heffter arrays and cycle decompositions of the complete graph.  We
have already discussed the connection between Heffter arrays and biembeddings of two cycle systems (cycle decompositions). In this section we will explicitly show how one can produce the two cycle decompositions $\mathcal{D}$ and $\mathcal{D'}$ and that these cycle decompositions have the property that if two edges are in the same cycle of $\mathcal{D}$, then they will be in different cycles in $\mathcal{D}'$.

First of all we recall some basic definitions about graphs and graph decompositions.
Good references about these concepts are \cite{W} and \cite{BEZ}, respectively.
Given a graph $\G$, by $V(\G)$ and $E(\G)$ we mean the vertex set and the edge set of $\G$, respectively,
and by $^\lambda \G$ the multigraph obtained from $\G$ by repeating each edge $\lambda$ times.
We denote by $K_{q \times r}$ the complete multipartite graph with $q$ parts each of size $r$,
obviously $K_{q \times 1}$ is nothing but the complete graph $K_q$. Finally, by $P_k$  we mean the path
of length $k$.

Let $\Gamma$ be a subgraph of a graph $K$. A $\Gamma$-\emph{decomposition} of $K$ is a set $\mathcal{D}$ of subgraphs
of $K$ isomorphic to $\G$, called \emph{blocks}, whose edges partition $E(K)$.
As already remarked in the previous section, a $C_k$-decomposition of $K_v$ is also called a $k$-\emph{cycle system of order} $v$.

\begin{definition}
Two graph decompositions $\mathcal{D}$ and $\mathcal{D}'$ of the same graph $K$ are  \emph{orthogonal}
if for every $B \in \mathcal{D}$ and every $B' \in \mathcal{D}'$, $B$ intersects $B'$ in at most one edge.
\end{definition}

We note the following obvious fact connecting the notions of biembedded cycle systems to orthogonal cycle systems.  From this fact, many of our earlier results on biembeddings of cycle systems  give results on orthogonal cycle systems.

\begin{proposition}\label{biembed.implies.orthog}
  If there is a simple biembedding
of a cycle system $ \mathcal{D}$ and a cycle system $ \mathcal{D}'$ in $K_v$, then $ \mathcal{D}$ and $ \mathcal{D'}$ are orthogonal cycle systems of order $v$.
\end{proposition}

If the vertices of $\G$ belong to a group $G$,
given $g \in G$, by $\G +g$ one means the graph whose vertex set is $V(\G)+g$ and whose edge set is $\{\{x+g,y+g\}\mid \{x,y\}\in E(\G)\}$.
An \emph{automorphism group} of a $\G$-decomposition $\mathcal{D}$ of $K$ is a group of bijections on $V(K)$ leaving $\mathcal{D}$ invariant.
A $\G$-decomposition of $K$ is said to be \emph{regular under a group} $G$ or $G$-\emph{regular} if it admits $G$
as an automorphism group acting sharply transitively on $V(K)$.
Heffter arrays are related to \emph{cyclic} decompositions, namely decompositions regular under a cyclic group.
We recall the following result.

\begin{proposition}
Let $K$ be a graph with $v$ vertices. A $\G$-decomposition $\mathcal{D}$ of $K$ is cyclic if and only if,
up to isomorphisms, the following conditions hold:
\begin{itemize}
\item[$(\rm{1})$ ] the set of vertices of $K$ is $\Z_v$;
\item[$(\rm{2})$ ] for all $B\in \mathcal{D}$,  $B+1\in \mathcal{D}$.
\end{itemize}
\end{proposition}

Clearly, to describe a cyclic decomposition it is sufficient to exhibit a complete system $\mathcal{B}$
of representatives for the orbits of $\mathcal{D}$ under the action of $\Z_v$. The elements of $\mathcal{B}$
are called \emph{base blocks} of $\mathcal{D}$.

One of most efficient tools applied for constructing regular decompositions is the \emph{difference method}, see \cite{AB}.
In particular, difference families over graphs have been introduced by Buratti \cite{B}, see also \cite{BP}.
Given a graph   $\G$ with
vertices in an additive group $G$, the \emph{list of differences} of $\G$
is the multiset $\Delta\G$ of all possible differences $x-y$ between two adjacent vertices of $\G$.
So, denoting by $D(\G)$ the set of \emph{di-edges} of $\G$, i.e., the set of all \underline{ordered}
pairs $(x,y)$ with $\{x,y\} \in E(\G)$, one can write:
$$\Delta \G = \{x-y \mid  (x,y)\in D(\Gamma)\}.$$

More generally, if $\mathcal{F}=\{\G_1,\G_2,\ldots,\G_{\ell}\}$ is a collection of graphs with vertices
in $G$, one sets $\Delta \mathcal{F}=\Delta \G_1 \cup \Delta \G_2 \cup \ldots \cup \Delta \G_{\ell}$, where in  the union
each element has to be counted with its multiplicity.

\begin{definition}
Let $J$ be a subgroup of a group $G$ and let $\G$ be a graph.
A collection $\mathcal{F}$ of graphs isomorphic to $\G$ and with vertices in $G$
is said to be a $(G,J,\G,\lambda)$-\emph{difference family} (briefly DF) \emph{over $G$ relative to $J$} if each element of
$G\setminus J$ appears exactly $\lambda$ times in $\Delta \mathcal{F}$ while no element of $J$ appears there.
\end{definition}

In the case of $J=\{0\}$ one simply says that $\mathcal{F}$ is a $(G,\G,\lambda)$-DF.
Also, if $t$ is a divisor of $v$, by writing $(v,t,\G,\lambda)$-DF one means a $(\Z_v,\frac{v}{t}\Z_v,\G,\lambda)$-DF
where $\frac{v}{t}\Z_v$ denotes the subgroup of $\Z_v$ of order $t$.
Analogously, $(v,\G,\lambda)$-DF means $(\Z_v,\G,\lambda)$-DF.
The connection between relative difference families and decompositions of a complete multipartite multigraph is given by the following result.

\begin{proposition}{\rm \cite{BP}}\label{DF-GD}
If $\mathcal{F}=\{\G_1,\G_2,\ldots,\G_{\ell}\}$ is a $(G,J,\G,\lambda)$-DF,
then $\mathcal{D}=\{\G_i+g \mid i=1,\ldots, \ell; g \in G\}$
is a $G$-regular $\G$-decomposition of $^\lambda K_{q \times r}$, where $q=|G:J|$ and $r=|J|$.
In particular, if $\mathcal{F}$ is $(v,\G,1)$-DF then $\mathcal{D}$ is a cyclic $\G$-decomposition of $K_{v}$.
\end{proposition}

Now we show how starting from a simple Heffter array one can construct a pair of orthogonal cyclic cycle decompositions of the complete graph $K_v$
or of
the cocktail party graph $K_{2\ell}-I$, namely the complete graph $K_{2\ell}$ minus the $1$-factor $I$ whose edges are
$[0,\ell],[1,\ell+1],[2,\ell+2],\ldots,[\ell-1,2\ell-1]$.
It is well known that the problem of finding necessary and sufficient conditions for cyclic cycle decompositions of $K_v$ and $K_{2\ell}-I$ has attracted much attention
(see, for instance, \cite{BDF,V,WF} and \cite{BGL,BR,JM2008,JM2017}, respectively).

First of all, from a simple Heffter array one can construct two cyclic difference families whose blocks are cycles,
as follows.
Let $H$ be a simple $H(m,n;h,k)$, hence every row and every column of $H$ admits a simple ordering.
Starting from a row of $H$ one can construct a cycle of length $h$ whose vertices are
the partial sums of the simple ordering of such a row. Denoted by $\Gamma_{i}$ the cycle
obtained from the $i$-th row $R_i$ of $H$, it is immediate to see that $\Delta\G_i =\pm R_i$.
Now, let $\mathcal{F}_R$ be the set of the $m$ cycles of length $h$ constructed from the $m$ rows of $H$ with this procedure.
By Condition (\rm{b}) of Definition \ref{def:Heffter}, we have $\Delta{\mathcal{F}_R}=\Z_{2mh+1}\setminus\{0\}$.
Hence $\mathcal{F}_R$ is a $(2mh+1, C_h,1)$-DF. Clearly, with the same technique, one can construct $n$ cycles of length $k$
starting from the columns of $H$. Let $\G'_j$ be the cycle obtained from the $j$-th column and let $\mathcal{F}_C$ be the set of these cycles.
Again by Condition (\rm{b}) of Definition \ref{def:Heffter}, we have $\Delta{\mathcal{F}_C}=\Z_{2nk+1}\setminus\{0\}$.
Hence $\mathcal{F}_C$ is a $(2nk+1, C_k,1)$-DF.

\begin{remark}\label{rem:ortho}
For every $\G_i \in \mathcal{F}_R$ and every $\G'_j \in \mathcal{F}_C$,
we see that $\G_i$  and $\G'_j$ contain at most one edge in common and hence  $|\Delta \G_i \cap \Delta \G'_j|\in\{0,2\}$. In particular, the intersection is empty if
the cell on the $i$-th row and the $j$-th column of $H$ is empty, while the intersection has size two otherwise.
\end{remark}

By Proposition \ref{DF-GD}, the cycles of $\mathcal{F}_R$ are the base blocks of a cyclic
$C_h$-decomposition $\mathcal{D}_R$ of $K_{2mh+1}$, while the cycles of $\mathcal{F}_C$ are the base blocks of a cyclic
$C_k$-decomposition $\mathcal{D}_C$ of $K_{2nk+1}$.
By Remark \ref{rem:ortho}, the
 decompositions $\mathcal{D}_R$ and $\mathcal{D}_C$ are orthogonal.

\begin{example}
Let $H$ be the $H(8;7)$ of {\rm Example \ref{ex:H87}} and let $\omega_i$'s and $\nu_i$'s be
 the simple orderings  given in {\rm Example \ref{ex:simple}}.
Considering the partial sums of the $\omega_i$'s ($\nu_i$'s, respectively) in $\Z_{113}$ we obtain the cycles $\G_i$'s ($\G'_i$'s, respectively):
\begin{center}
\begin{footnotesize}
$\begin{array}{lcl}
\G_1= ( 8, 33, 49, 22, -7, 24, 0); & \quad &\G'_1= (8,47,30,-8,-51,-9,0);\\
\G_2= (-17, -23, -51, -28, -2, 30, 0); & \quad &\G'_2= ( 16,10,-35,-45,3,52,0 );\\
\G_3= (39, 29, 24, 57, 72, 37, 0); & \quad &\G'_3= ( 23,18,65,47,1,-50,0);\\
\G_4= (-38, -56, -49, -85, -74, -40, 0); & \quad &\G'_4= (25,-3,12,19,-34,-56,0);\\
\G_5= (-43, -88, -41, -63, -60, -19,  0); & \quad &\G'_5= (-27,-1,10,13,68,54,0);\\
\G_6= ( 42, 90, 44, 30, 32, -12, 0); & \quad &\G'_6= ( -21,11,44,8,27,29,0);\\
\G_7= (20, -31, -84, -29, -50, -49, 0);& \quad  &\G'_7= (-13,-43,-78,-44,-32,-31,0);\\
\G_8= (-52, -43, 7, 63, 9, -4, 0); & \quad &\G'_8= (-37,-61,-21,20,-24,-4,0).
\end{array}$
\end{footnotesize}

\end{center}
Then $\mathcal{F}_R=\{\G_1,\ldots,\G_8\}$ and $\mathcal{F}_C=\{\G'_1,\ldots,\G'_8\}$ are two sets of base cycles
of a pair of orthogonal cyclic $C_7$-decompositions of  $K_{113}$.
\end{example}

From the arguments in this section we get the following results.  Note the similarity of Proposition \ref{orth.systems} to Theorem \ref{biembed2}, however we see that for the application of Heffter arrays to orthogonal cycle decompositions it is not required that the orderings of the rows and columns of the Heffter array be compatible.  For results on orthogonal cycle decompositions see \cite{BCP} and the references therein.

\begin{proposition}
If there exists a simple $H(m,n;h,k)$, then there exist
a $(2mh+1,C_h,1)$-DF and a $(2nk+1, C_k,1)$-DF.
\end{proposition}

\begin{proposition}\label{orth.systems}
If there exists a simple $H(m,n;h,k)$, then there exist
 a cyclic $C_h$-decomposition and
a cyclic $C_k$-decomposition both of $K_{2nk+1}$,
and these decompositions  are orthogonal.
\end{proposition}

Obviously if the array is globally simple one can immediately construct the cycles starting from the rows and from the columns,
without providing the simple orderings.
Also note that the elements of an $H(m,n;h,k)$ belong to $\Z_{2nk+1}$, hence when we say that an ordering of a row or
of a column is simple, we mean that its partial sums are pairwise distinct modulo $2nk+1$.
In \cite{CMPPHeffter}, the authors showed that it is useful to find orderings whose partial sums are pairwise distinct both
modulo $2nk+1$ and modulo $2nk+2$. To this end they give the following definition.

\begin{definition}An $SH^*(m,n;h,k)$ is a globally simple $H(m,n;h,k)$
such that the natural ordering of each row and each column is simple also modulo $2nk+2$.
\end{definition}
The arrays of Theorem \ref{thm:SH} have this additional property.
The usefulness of these arrays is explained by the following proposition.

\begin{proposition}\label{pr:ccp}{\rm \cite{CMPPHeffter}}
If there exists an $SH^*(m,n;h,k)$, then there exist:
\begin{itemize}
\item[{\rm (1)}] a cyclic $C_h$-decomposition  and a cyclic $C_k$-decomposition both of $K_{2nk+1}$ and these decompositions are orthogonal;
\item[{\rm (2)}] a cyclic $C_h$-decomposition  and a cyclic $C_k$-decomposition both of $K_{2nk+2}-I$ and these decompositions are orthogonal.
\end{itemize}
\end{proposition}

\begin{example}\label{10.8}
The following is an $SH^*(10;8)$ over $\Z_{161}$:
\begin{center}
\begin{footnotesize}
$\begin{array}{|r|r|r|r|r|r|r|r|r|r|}\hline
77 & 80 & -78 & -71 & -70 & -79 &  &  & 69 & 72 \\\hline
 &  & -17 & -20 & -25 & -28 & 26 & 19 & 18 & 27 \\\hline
5 & 8 & 13 & 16 & -14 & -7 & -6 & -15 &  &  \\\hline
34 & 43 &  &  & -33 & -36 & -41 & -44 & 42 & 35 \\\hline
 &  & 21 & 24 & 29 & 32 & -30 & -23 & -22 & -31 \\\hline
58 & 51 & 50 & 59 &  &  & -49 & -52 & -57 & -60 \\\hline
-38 & -47 &  &  & 37 & 40 & 45 & 48 & -46 & -39 \\\hline
-73 & -76 & 74 & 67 & 66 & 75 &  &  & -65 & -68 \\\hline
-62 & -55 & -54 & -63 &  &  & 53 & 56 & 61 & 64 \\\hline
-1 & -4 & -9 & -12 & 10 & 3 & 2 & 11 &  &  \\\hline
\end{array}$
\end{footnotesize}
\end{center}
By the partial sums in $\Z_{161}$ of the natural  orderings of the rows (columns, respectively)  we obtain the cycles
$\G_i$'s ($\G'_i$'s, respectively):
\begin{center}
\begin{footnotesize}
$\begin{array}{lcl}
\G_1= ( 77,   157,    79,     8,   -62,  -141,   -72,     0 ); & \;\; &\G'_1= (77,  82,  116, 13,   136,    63,     1,     0) ;\\
\G_2= (-17,   -37,   -62,   -90,   -64,   -45,   -27,     0 ); & \;\; &\G'_2= (80,    88,131,21, 135,    59,     4,     0  );\\
\G_3= (5,    13,    26,    42,    28,    21,    15,     0 ); & \;\; &\G'_3= ( -78,   -95,   -82,   -61,   -11,    63,   9,     0 ); \\
\G_4= (34,    77,    44,     8,   -33,   -77,   -35,     0); & \;\; &\G'_4= (-71,   -91,   -75,   -51,     8,    75,    12,   0);\\
\G_5= (21,    45,    74,   106,    76,    53,    31,0); & \;\;&\G'_5= (-70,   -95,  -109,  -142,  -113,    -76, -10,     0 );\\
\G_6=(58,   109,   159,   57,   8,   117,    60,     0); & \;\; &\G'_6= (-79,  -107,  -114,  -150,  -118,   -78, -3,     0  );\\
\G_7= (-38,   -85,   -48,    -8,    37,    85,    39,     0 );& \;\;  &\G'_7= ( 26,    20,   -21,   -51,  -100,   -55,  -2,     0 );\\
\G_8= ( -73,  -149,   -75,    -8,    58,   133,    68,     0 ); & \;\; &\G'_8= (19,     4,   -40,   -63,  -115,   -67, -11,     0 );\\
\G_9= (-62,  -117,  -10,  -73,  -20,  -125,   -64,     0 ); &\;\;&\G'_9=( 69,    87,   129,   107,    50,     4,   -61,     0 );\\
\G_{10}=(-1,    -5,   -14,   -26,   -16,   -13,   -11,     0); & \;\; &\G'_{10}=(72,    99,   134,   103,    43,     4,   -64,     0 ).
\end{array}$
\end{footnotesize}
\end{center}
Then $\mathcal{F}_R=\{\G_1,\ldots,\G_{10}\}$ and $\mathcal{F}_C=\{\G'_1,\ldots,\G'_{10}\}$ are two sets of base cycles
of a pair of orthogonal cyclic $C_8$-decompositions of  $K_{161}$.\\
Analogously, if we consider the partial sums of each row (column, respectively) in $\Z_{162}$, we obtain the cycles
$\widetilde{\G}_i$'s ($\widetilde{\G}'_j$'s, respectively) where:
\begin{center}
\begin{footnotesize}
$\begin{array}{lcl}
\widetilde{\G}_i= \G_i, \quad i\ne 6,9;  &\quad&  \widetilde{\G}'_j=\G_j, \quad j\ne 1,2;\\
\widetilde{\G}_6= (58,   109,   159,   56,   7,   117,    60,     0);& \quad &\widetilde{\G}'_1=
(77,    82,   116,   12,   136,    63,     1,     0);\\
\widetilde{\G}_9=(-62,  -117,  -9,  -72,  -19,  -125,   -64,     0 ); & \quad &\widetilde{\G}'_2= (80,    88,   131,   20,   135,    59,     4,     0  ).
\end{array}$
\end{footnotesize}
\end{center}
 Now $\widetilde{\mathcal{F}_R}=\{\widetilde{\G}_1,\ldots,\widetilde{\G}_{10}\}$ and
$\widetilde{\mathcal{F}_C}=\{\widetilde{\G}'_1,\ldots,\widetilde{\G}'_{10}\}$ are two sets of base cycles
of a pair of orthogonal cyclic $C_8$-decompositions of the cocktail party graph $K_{162}-I$.
\end{example}

\section{Variants and generalizations}\label{sec:variants}

In  this section we will discuss
several of the variants and generalizations of Heffter arrays that have been researched in earlier publications.

\subsection{Weak Heffter arrays}

In the original paper on Heffter arrays \cite{A}, Archdeacon also proposed the following weaker type of Heffter array.  These weaker Heffter arrays  again have applications  to biembedding cycle systems as well as to orthogonal cycle systems.
Two Heffter systems $D_h=D(2mh+1,h)$ and $D_k=D(2nk+1,k)$ with $mh = nk$ are {\em weakly sub-orthogonal}
if the $i$-th part of $D_h$ has at most one element $a_{i,j}$ such that either $a_{i,j}$ or $-a_{i,j}$
is in the $j$-th part of $D_k$. Form a {\em weak Heffter array} $H(m,n;h,k)$ by placing $a_{i,j}$ in
the cell on the $i$-th row and the $j$-th column.
The upper sign on $\pm$ or $\mp$ is the {\em row sign} corresponding to its sign on $a_{i,j}$ in $D_h$,
the lower sign is the {\em column sign} used in $D_k$.
Using the row signs we get row sums 0 and the column signs give column sums 0.

The following is a formal definition of a weak Heffter array.
\begin{definition}\label{def:weak}
 A \emph{weak Heffter array} $H(m,n;h,k)$ is an $m \times n$ matrix $A$  such that:
\begin{itemize}
  \item[$(\rm{a_1})$] each row contains $h$ filled cells and each column contains $k$ filled cells;
\item[$(\rm{b_1})$]  for every $x \in \mathbb{Z}_{2nk+1} \setminus \{0\}$, there is exactly one cell of $A$ whose element is one of the following:
$x,-x,\pm x,\mp x$, where the upper sign on $\pm$ or $\mp$ is the row sign  and
the lower sign is the column sign;
\item[$(\rm{c_1})$] the elements in every row and column (with the corresponding sign) sum  to $0$ in $\mathbb{Z}_{2nk+1}$.
\end{itemize}
\end{definition}

\begin{example}
Consider the following Heffter systems on a half-set of $\Z_{25}$:
$$D_4=\{\{1,-7-6,12\},\{2,-4,10,-8\},\{-3,-11,9,5\}\}$$
and
$$D_3=\{\{1,2,-3\},\{-7,-4,11\},\{-6,-10,-9\},\{12,8,5\}\}.$$
Note that they are weakly sub-orthogonal and starting from these two systems one gets
the following  $3 \times 4$ weak Heffter array over $\Z_{25}$ without empty cells :
$$
\begin{array}{|r|r|r|r|} \hline
1 & -7 & -6 & 12 \\ \hline
2 & -4 & \pm 10 & \mp 8 \\ \hline
-3 & \mp 11 & \pm 9 & 5 \\ \hline
\end{array}
$$
\end{example}

To date there are no published papers specifically on weak Heffter arrays, but one is in preparation, see \cite{weak}.
In \cite{A} Archdeacon showed how these arrays are related to current graphs and hence to biembeddings
(in this case on a {\em non-orientable} surface).
Also, by reasoning analogous to that in Section \ref{sec:DF}, it is easy to see that the following holds.

\begin{proposition}\label{weak}
If there exists a simple weak $H(m,n;h,k)$, then there exist
a $(2mh+1,C_h,1)$-DF and a $(2nk+1, C_k,1)$-DF.
\end{proposition}

\begin{proposition}
If there exists a simple weak $H(m,n;h,k)$, then there exist
a cyclic $C_h$-decomposition  and
 a cyclic $C_k$-decomposition both of $K_{2nk+1}$.
These decompositions  are orthogonal.
\end{proposition}

\subsection{$\lambda$-fold relative Heffter arrays}

In \cite{RelH,CPEJC} the authors proposed a natural generalization of Heffter arrays,
which is related to signed magic arrays, difference families,
graph decompositions and biembeddings.  In this generalization, the array is missing the elements
of a subgroup $J\subset \Z_v$ and each element of $ \Z_v \setminus J$ occurs exactly $\lambda$ times in the array.  We give the following definition.

\begin{definition}\label{def:lambdaRelative}
Let $v=\frac{2nk}{\lambda}+t$ be a positive integer,
where $t$ divides $\frac{2nk}{\lambda}$,  and
let $J$ be the subgroup of $\Z_{v}$ of order $t$.
 A  $\lambda$-\emph{fold Heffter array $A$ over $\Z_{v}$ relative to $J$}, denoted by $^\lambda  H_t(m,n; h,k)$, is an $m\times n$ array
 with elements in $\Z_{v}\setminus J$ such that:
\begin{itemize}
\item[$(\rm{a_2})$] each row contains $h$ filled cells and each column contains $k$ filled cells;
\item[$(\rm{b_2})$] the multiset $\{\pm x \mid x \in A\}$ contains  each element of $\Z_v\setminus J$ exactly $\lambda$ times;
\item[$(\rm{c_2})$] the elements in every row and column sum to $0$ in $\Z_v$.
\end{itemize}
\end{definition}

The terminology is due to their connection with \emph{relative} difference families.
Trivial necessary conditions for the existence of a $^\lambda H_t(m,n; h,k)$ are $mh=nk$, $2\leq h \leq n$ and $2\leq k \leq m$.
It is easy to see that condition $(\rm{b_2})$ of Definition \ref{def:lambdaRelative} asks that
for every $x\in \Z_{v}\setminus J$ with $x$ different from the involution the number of  occurrences of $x$ and $-x$ in the array is $\lambda$,
while if the involution exists and it does not belong to $J$, then it has to appear exactly $\frac{\lambda}{2}$ times,
also no element of $J$ appears in the array.

If $^\lambda H_t(m,n; h,k)$ is a square array, then it is denoted by $^\lambda H_t(n;k)$.
If $\lambda=t=1$ one finds again the classical concept of  Heffter array, i.e. an
$^1 H_1(m,n; h,k)$ is just an $H(m,n; h,k)$.
In general if $\lambda=1$, then it is omitted; the same holds for $t$.
 Note also that if $\lambda=1$, then $3\leq h \leq n$ and $3\leq k \leq m$.

\begin{example}\label{ex:1}
A $^3H_2(4;3)$ whose elements belong to $\Z_{10}$ and a $^4H_4(4;2)$ whose elements belong to $\Z_8$:
  \center{
  $\begin{array}{|r|r|r|r|}
\hline  & 1 & 2  & -3  \\
\hline 4  &  & 4 & 2  \\
\hline -3  & 2 &  &  1 \\
\hline   -1 & -3 & 4 &  \\
\hline
\end{array}$
\quad\quad\quad\quad
  $\begin{array}{|r|r|r|r|}
\hline  1 & -1 &   &   \\
\hline -1  & 1 &  &   \\
\hline   &  & 3 &  -3 \\
\hline   &  & -3 & 3 \\
\hline
\end{array}$
}
\end{example}
As was the case for Heffter arrays, a $\lambda$-fold relative Heffter array is called \emph{integer} if condition ($\rm{c_2}$) in Definition \ref{def:lambdaRelative}
is strengthened so that the elements in every row and every column, seen as integer in
$\pm\left\{1,2,\ldots, \lfloor\frac{v}{2}\rfloor\right\}$, sum to $0$ in $\Z$.
The $^3H_2(4;3)$ of Example \ref{ex:1} is not an integer Heffter array, while the $^4H_4(4;2)$ in the same example is integer.

\begin{example}\label{ex:inv}
 An integer $^2H(5;3)$ over $\Z_{16}$
(note that the involution
 appears exactly once in the array) and an integer $^2 H(6;4)$ over $\Z_{25}$:

\center{
  $\begin{array}{|r|r|r|r|r|}
\hline 3 & 5 & -8  & &  \\
\hline & 2 & 5 & -7 &   \\
\hline & & 3  & 1 &  -4 \\
\hline  -4 &  &  & 6 & -2 \\
\hline   1 & -7 & &  & 6 \\
\hline
\end{array}$ \hspace{1 cm}
  $\begin{array}{|r|r|r|r|r|r|}
\hline & & 1 & -2  & -5 & 6  \\
\hline & & -3  & 4 & 7 & -8  \\
\hline -9 & 10  &  &  &  1 & -2 \\
\hline   11 & -12 & & &  -3 & 4   \\
\hline 5 & -6 &-9 & 10  &  &   \\
\hline  -7 & 8 & 11 & -12 & &   \\
\hline
\end{array}$
 } 

\end{example}

\subsubsection{Necessary conditions}

In \cite{CPEJC} the authors established some necessary conditions for the existence of a $^\lambda H_t(m,n;h,k)$.
Firstly, it is trivial to note that Conditions ($\rm{a_2})$ and ($\rm{b_2})$ of Definition \ref{def:lambdaRelative} imply that
if there exists a $^{\lambda} H_t(m,n;h,k)$ with either $h=2$ or $k=2$, then $\lambda$ has to be even.

\begin{proposition}\label{prop:necctrivial}{\rm \cite{CPEJC}}
Suppose that there exists a $^\lambda H_t(m,n;h,k)$ and set $v=\frac{2nk}{\lambda}+t$.
 If either $v$ is odd or $v$ and $t$ are even,  then $\lambda$ has to be a divisor of $nk$.
\end{proposition}

\begin{proposition}\label{prop:nonexistence}{\rm \cite{CPEJC}}
If $\lambda\equiv 2 \pmod 4$, $v=\frac{2nk}{\lambda}+t\equiv2\pmod 4$ and $t$ is odd, then a  $^{\lambda} H_t(m,n;h,k)$
cannot exist.
\end{proposition}

About the integer case there are the following more restrictive conditions.
\begin{proposition}\label{prop:necc}{\rm \cite{RelH,CPEJC}}
Suppose that there exists an integer $^\lambda H_t(m,n;h,k)$ with $\lambda$ odd.
\begin{itemize}
\item[$(\rm{1})$] If $t$ divides $\frac{nk}{\lambda}$, then
$\frac{nk}{\lambda}\equiv 0 \pmod 4$ or $\frac{nk}{\lambda}\equiv -t \equiv \pm 1\pmod 4.$
\item[$(\rm{2})$] If $t=\frac{2nk}{\lambda}$, then $h$ and $k$ must be even.
\item[$(\rm{3})$] If $t\neq \frac{2nk}{\lambda}$ does not divide $\frac{nk}{\lambda}$, then
$\frac{2nk}{\lambda}+t\equiv 0 \pmod 8.$
\end{itemize}
\end{proposition}

It is known that these conditions are not sufficient. For instance in \cite{RelH} it is proved that
there is no integer $H_{3n}(n;3)$ for $n\geq 3$ and no integer $H_8(4;3)$, moreover in \cite{weak}
it is shown the nonexistence of an integer $H_t(4;3)$ for $t=4,6$.

\subsubsection{Integer relative Heffter arrays with $\lambda=1$}
\vspace{1 em}
Turning to the existence question now, the first class investigated is that of  integer relative Heffter arrays with $\lambda=1$.

\begin{theorem}{\rm \cite{RelH}}
Let $3 \leq k \leq n$ with $k \neq 5$. There exists an integer $H_k (n; k)$ if and only if
one of the following holds:
\begin{itemize}
  \item [$(\rm{1})$] $k$ is odd and $n \equiv 0, 3 \pmod 4$;
  \item [$(\rm{2})$] $k \equiv 2 \pmod 4$ and $n$ is even;
  \item [$(\rm{3})$] $k \equiv 0 \pmod 4$.
\end{itemize}
Furthermore, there exists an integer $H_5(n; 5)$ if $n \equiv 3 \pmod 4$ and it does not exist if
$n \equiv 1, 2 \pmod 4$.
\end{theorem}
The case $k=5$ and $n\equiv0\pmod 4$ is still open. For these values of the parameters integer
weak Heffter arrays
have been constructed in \cite{weak}.

\begin{theorem} {\rm \cite{CPPBiembeddings}}
There exists an integer globally simple $H_t(n;k)$ in the following cases:
\begin{itemize}
  \item [$(\rm{1})$] $k=3$, $t=n,2n$ for all odd $n$;
  \item [$(\rm{2})$] $t=k=3,5,7,9$ for all $n\equiv 3 \pmod4$.
\end{itemize}
\end{theorem}

\begin{theorem} {\rm \cite{MP}}
Let $m, n, h, k$ be such that $4 \leq h \leq n$, $4 \leq k \leq m$ and $mh = nk$. Let $t$
be a divisor of $2nk$.
\begin{itemize}
  \item [$(\rm{1})$] If $h, k \equiv 0 \pmod 4$, then there exists an integer $H_t(m, n; h, k)$.
  \item [$(\rm{2})$] If $h \equiv 2 \pmod 4$ and $k \equiv 0 \pmod 4$, then there exists an integer $H_t(m, n; h, k)$ if and
only if $m$ is even.
  \item [$(\rm{3})$] If $h \equiv 0 \pmod 4$ and $k \equiv 2 \pmod 4$, then there exists an integer $H_t(m, n; h, k)$ if and
only if $n$ is even.
  \item [$(\rm{4})$] Suppose that $m$ and $n$ are both even. If $h, k \equiv 2 \pmod 4$, then there exists an integer
$H_t(m, n; h, k)$.
\end{itemize}
\end{theorem}

\begin{theorem} {\rm \cite{MP3}}
Let $m, n, h, k$ be such that $3 \leq h \leq n$, $3 \leq k \leq m$ and $mh = nk$.
Suppose $d = \gcd(h, k) \geq 3$.
\begin{itemize}
  \item[$(\rm{1})$] If $d \equiv 1 \pmod 4$ and $nk \equiv 3 \pmod 4$, then there exists an integer $H_d(m, n; h, k)$.
  \item[$(\rm{2})$] If $d \equiv 3 \pmod 4$ and $nk \equiv 0, 1 \pmod 4$, then there exists an integer $H_d(m, n; h, k)$.
  \item[$(\rm{3})$]  If $d = 3$ and $nk$ is odd, then there exists an integer $H_\frac{nk}{3}(m, n; h, k)$ and an integer
$H_\frac{2nk}{3}(m, n; h, k)$.
\end{itemize}
\end{theorem}

\subsubsection{  $\lambda$-fold relative Heffter arrays}
The following are the known existence results on $\lambda$-fold relative Heffter arrays.

\begin{theorem} {\rm \cite{MP3}}
Let $m, n, h, k$ be such that $3 \leq h \leq n$, $3 \leq k \leq m$ and $mh = nk$.
Suppose $d = \gcd(h, k) \geq 3$.
\begin{itemize}
  \item[$(\rm{1})$] If $d = 3$ and $nk \equiv 3 \pmod 4$, then there exists a  $^2 H(m, n; h, k)$.
  \item[$(\rm{2})$] If $d = 3$ and $nk \equiv 1 \pmod 4$, then there exists a  $^3 H(m, n; h, k).$
  \item[$(\rm{3})$]  If $d = 3$, $nk$ is odd, and $\lambda$ divides $n$, then there exists a  $^\lambda H_\frac{n}{\lambda}(m, n; h, k).$
  \item[$(\rm{4})$] If $d = 3$, $nk$ is odd, and $\lambda$ divides $2n$, then there exists a  $^\lambda H_\frac{2n}{\lambda}(m, n; h, k).$
  \item[$(\rm{5})$]  If $d = 5$ and $nk \equiv 3 \pmod 4$, then there exists a  $^2 H(m, n; h, k).$
\end{itemize}
\end{theorem}

\begin{theorem} {\rm \cite{MP2}}
Let $m, n, h, k$ be  such that $4 \leq h \leq n$, $4 \leq  k \leq m$ and $mh = nk$. Let
$\lambda$ be a divisor of $2nk$ and let $t$ be a divisor of $\frac{2nk}{\lambda}$.
 There exists an integer $^\lambda H_t(m, n; h, k)$ in each of the following cases:
 \begin{itemize}
  \item[$(\rm{1})$] $h, k \equiv 0 \pmod 4$;
 \item[$(\rm{2})$]  $h \equiv 2 \pmod 4$ and $k \equiv 0 \pmod 4$;
 \item[$(\rm{3})$] $h \equiv 0 \pmod 4$ and $k \equiv 2 \pmod 4$;
 \item[$(\rm{4})$] $h, k \equiv 2 \pmod 4$ and $m, n$ both even.
  \end{itemize}
\end{theorem}

\begin{theorem}\label{lambda2}{\rm \cite{CPEJC}}
There exists a $^2 H(n;k)$ if and only if $n\geq k\geq 3$ with $nk\not\equiv 1 \pmod 4$.
\end{theorem}

\begin{theorem} \label{thm:noemptycell}{\rm \cite{CPEJC}}
There exists a $^2 H(m,n;n,m)$ if and only if $mn\not\equiv 1 \pmod 4$ and $m,n$ are not both equal to $2$.
\end{theorem}

\begin{theorem}\label{thm:skeven}{\rm \cite{CPEJC}}
Let $h,k$ be even positive integers.
There exists a $^2 H(m,n;h,k)$ if and only if $mh=nk$, $2\leq h\leq n$, $2\leq k \leq m$ and $h,k$ are not both equal to $2$.
Moreover these arrays are integer when:
\begin{itemize}
\item[$(\rm{1})$] $h,k\neq 2$;
\item[$(\rm{2})$] $h=2$ and either $n=2$ and $m=k\equiv 0\pmod 4$ or $n,k\geq3$;
\item[$(\rm{3})$] $k=2$ and either $m=2$ and $n=h\equiv 0\pmod 4$ or $m,h\geq 3$.
\end{itemize}
\end{theorem}

\begin{proposition}\label{lambdak}{\rm \cite{CPEJC}}
Let $n,k,\lambda$ be positive integers with $n\geq k\geq 3$ and $\lambda$ divisor of $k$.
  Then there exists a $^\lambda H_{\frac{k}{\lambda}}(n;k)$ if one of the following is satisfied:
  \begin{itemize}
  \item[$(\rm{1})$] $k=5$ and $n\equiv 3 \pmod 4$;
  \item[$(\rm{2})$] $k\neq 5$ odd and $n\equiv 0,3 \pmod 4$.
  \end{itemize}
\end{proposition}

\begin{proposition}\label{prop:n}{\rm \cite{CPEJC}}
  For any odd integer $n\geq 3$ and for any divisor $\lambda$ of $n$ there exists a $^\lambda H_\frac{n}{\lambda}(n;3)$.
\end{proposition}

\begin{proposition}\label{prop:2n}{\rm \cite{CPEJC}}
  For any odd integer $n\geq 3$ and for any divisor $\lambda$ of $2n$ there exists a $^\lambda H_\frac{2n}{\lambda}(n;3)$.
\end{proposition}

\begin{proposition}{\rm \cite{CPEJC}}
Let $m,n,h,k,t,\lambda$ be  such that $mh=nk$, $4\leq h\leq n$, $4\leq k\leq m$, $t$ is a divisor of $2nk$ and
$\lambda$ is a divisor of $t$. A $^\lambda H_{\frac{t}{\lambda}}(m,n;h,k)$ exists in each of the following cases:
\begin{itemize}
  \item[$(\rm{1})$] $h\equiv k\equiv 0\pmod 4$;
  \item[$(\rm{2})$] either $h\equiv 2\pmod 4$ and $k\equiv 0\pmod 4$ or
 $h\equiv 0\pmod 4$ and $k\equiv 2\pmod 4$;
  \item[$(\rm{3})$] $h\equiv k \equiv 2\pmod 4$ and $m$ and $n$ are even.
\end{itemize}
\end{proposition}

Some of previous results have been obtained thanks to one of the following two theorems
which allow one to obtain new arrays from  existing ones.

\begin{theorem}\label{thm:proiezione}{\rm \cite{CPEJC}}
If there exists an $^\alpha H_t(m,n;h,k)$ then, for any divisor $\lambda$
of $t$, there exists a $^{\lambda\alpha} H_\frac{t}{\lambda}(m,n;h,k)$.
\end{theorem}

\begin{theorem}{\rm \cite{MP3}}
Let $m, n, h, k$ be such that $2 \leq h \leq n$, $2 \leq k \leq m$ and $mh = nk$.
Let $\lambda$ be a divisor of $2nk$ and let $t$ be a divisor of $\frac{2nk}{\lambda}$.
 Set $d = \gcd(h, k)$. If there exists a (integer) diagonal $^\lambda H_t(\frac{nk}{d};d)$,
 then there exists a (integer) $^\lambda H_t(m, n; h, k)$.
\end{theorem}

Other existence results can be found using the following recursive construction.
\begin{proposition}\label{prop:recursive}{\rm \cite{CPEJC}}
If there exists an (integer) $H_t(m,n;h,k)$ then there exists an (integer) $^{\alpha_1\lambda_1}H_t(\lambda_1m,\lambda_2n;\alpha_1h,\alpha_2k)$
for all positive integers $\alpha_1\leq \lambda_2$, $\alpha_2\leq\lambda_1$ such that $\alpha_1\lambda_1=\alpha_2\lambda_2$,
$\alpha_1h\leq \lambda_2 n$ and $\alpha_2k\leq \lambda_1 m$.
\end{proposition}

\subsubsection{Connection between relative Heffter arrays and biembeddings and graph decompositions}

In \cite{RelH,CPEJC} the authors investigated the connection of these generalized arrays with relative difference families,
graph decompositions and biembeddings.
To illustrate the results we need some definitions.
Given a partially filled array $A$ by \emph{skeleton} of $A$, denoted by $skel(A)$, one means the set of the filled
 positions of $A$, while by $\E(A)$ one denotes the multiset of the elements of $skel(A)$.
 Analogously, by $\E(R_i)$ and $\E(C_j)$ we mean the multisets of elements
 of the $i$-th row and of the $j$-th column, respectively, of $A$.
Now we have to define a simple ordering of a \emph{multiset}.
We consider  the elements of $\E(A)$  indexed by the set $skel(A)$ (analogously for the rows and the columns).
For example let $A$  be the $^4 H_4(4;2)$ constructed in Example \ref{ex:1}.
Here $\E(A)=\{-1,-1,1,1,-3,-3,3,3\}$, we can view this multiset as the set
$\{t_b \mid b\in skel(A)\}=\{t_{(1,1)},t_{(1,2)},t_{(2,1)},t_{(2,2)},t_{(3,3)},t_{(3,4)},t_{(4,3)},t_{(4,4)}\}$,
where $t_{(1,1)}=1$, $t_{(1,2)}=-1$, $t_{(2,1)}=-1$, $t_{(2,2)}=1$, $t_{(3,3)}=3$, $t_{(3,4)}=-3$,
$t_{(4,3)}=-3$ and $t_{(4,4)}=3$.
 Given a finite multiset $T=[t_{b_1},\dots,t_{b_k}]$ whose (not necessarily distinct) elements are indexed by a set $B=\{b_1,\dots,b_k\}$ and a cyclic permutation $\alpha$ of $B$ we say that the list $\omega_{\alpha}=(t_{\alpha(b_1)},t_{\alpha(b_2)},\ldots,t_{\alpha(b_k)})$ is the ordering of the elements of $T$ associated to $\alpha$. In case the elements of $T$ belong to an abelian group $G$, we define $s_i=\sum_{j=1}^i t_{\alpha(b_j)}$, for any $i\in\{1,\ldots,k\}$, to be the $i$-th partial sum of $\omega_{\alpha}$ and we set $\S(\omega_{\alpha})=(s_1,\ldots,s_k)$.
The ordering $\omega_{\alpha}$ is said to be \emph{simple} if $s_b\neq s_c$ for all $1\leq b <  c\leq k$.

\begin{proposition}
If there exists a simple $^\lambda H_t(m,n;h,k)$, then there exist
a $(2mh+t,t,C_h,\lambda)$-DF and a $(2nk+t,t, C_k,\lambda)$-DF.
\end{proposition}

\begin{proposition}
If there exists a simple  $^\lambda H_t(m,n;h,k)$,  then there exist a cyclic $C_h$-decomposition $\mathcal{D}$  and
a cyclic $C_k$-decomposition $\mathcal{D}'$ both of $^\lambda K_{\frac{2nk+t}{t}\times t}$.
If $\lambda=1$ the decompositions $\mathcal{D}$ and $\mathcal{D}'$ are orthogonal.
\end{proposition}

Also for these generalized arrays to get biembeddings one has to find compatible orderings.
Note that in Section \ref{sec:biemb} we have defined compatible orderings only when $\E(A)$ is a set.
Given an $m \times n$ partially filled  array $A$, by $\alpha_{R_i}$ we will denote a cyclic permutation of the nonempty cells of the $i$-th row
$R_i$ of $A$ and by $\omega_{R_i}$ (omitting the dependence on $\alpha_{R_i}$) the associated ordering of $\E(R_i)$. Similarly, for the $j$-th column $C_j$ of $A$, we define $\alpha_{C_j}$ and $\omega_{C_j}$.

 If for any $i\in\{1,\ldots, m\}$ and for any $j\in\{1,\ldots,n\}$, the orderings $\omega_{R_i}$ and $\omega_{C_j}$ are simple, we define the permutation $\alpha_r$ of $skel(A)$
 by $\alpha_{R_1}\circ \ldots \circ\alpha_{R_m}$ and we say that the associated action $\omega_r$ on $\E(A)$ is the simple ordering for the rows.
 Similarly we define $\alpha_c=\alpha_{C_1}\circ \ldots \circ\alpha_{C_n}$ and the simple ordering $\omega_c$ for the columns.
Now, given a  $\lambda$-fold relative Heffter array $^\lambda H_t(m,n; h,k)$, say $A$, then the orderings $\omega_r$ (associated to the permutation $\alpha_r$) and $\omega_c$ (associated to the permutation $\alpha_c$) are said to be
\emph{compatible} if $\alpha_c \circ \alpha_r$ is a cycle of length $|skel(A)|$.
In \cite{CPEJC} the authors showed that the necessary conditions of Theorem \ref{thm:compatible}
hold more in general for a  $^\lambda H_t(m,n; h,k)$.

\begin{theorem}\label{thm:biembedding}{\rm \cite{CPEJC}}
Let $A$ be a  $\lambda$-fold relative Heffter array $^\lambda H_t(m,n; h,k)$ that is simple with respect to the compatible orderings $\omega_r$ and $\omega_c$.
Then there exists a cellular biembedding of a cyclic  $C_h$-decomposition  and
a cyclic $C_k$-decomposition both of $^\lambda K_{(\frac{2nk}{\lambda t}+1)\times t}$ into an orientable surface.
\end{theorem}

\subsection{Signed Magic Arrays}

The notion of a signed magic array was introduced by Khodkar, Schulz and Wagner   \cite{KSW} in 2017.  These are intimately related to (relative) Heffter arrays.  We begin with the definition.

\begin{definition}\label{def:SMR}
A \emph{signed magic array} $SMA(m,n;h,k)$ is an $m\times n$ array with entries from $X$,
where $X=\left\{0,\pm1,\pm2,\ldots, \pm \frac{nk-1}{2}\right\}$ if $nk$ is odd and $X=\left\{\pm1,\pm2,\ldots,\pm \frac{nk}{2}\right\}$
if $nk$ is even, such that:
\begin{itemize}
\item[$(\rm{a_3})$] each row contains $h$ filled cells and each column contains $k$ filled cells;
\item[$(\rm{b_3})$] every integer from the set $X$ appears exactly once in the array;
\item[$(\rm{c_3})$] the sum of each row and of each column is $0$ in $\Z$.
\end{itemize}
\end{definition}

\begin{remark}\label{rem:magicHeffter}
As was shown in \cite{CPEJC}, a signed magic array $SMA(m,n;h,k)$ with $nk$ even is an integer $^2 H(m,n;h,k)$, while
in general the converse is not true. For instance, the $^2 H(6;4)$ of Example \ref{ex:inv} is not a signed magic array.
On the other hand,  in the particular case in which either $h=2$ or $k=2$ an
integer $^2 H(m,n;h,k)$ is precisely an $SMA(m,n;h,k)$, see \cite{CPEJC} for details.
\end{remark}

Existence results on signed magic arrays can be found in \cite{KE,KL20,KL,KLE,KSW,MP2,MP3}. We note that nearly all the results and techniques used on signed magic arrays are analogous to those for Heffter arrays.
In particular, the existence of an $SMA(m, n; h, k)$ has been determined in the square case and when the array has no empty cells, as shown in the following two theorems.

\begin{theorem}{\rm \cite{KSW}}
There exists an $SMA(n, n; k, k)$ if and only if either $n = k = 1$ or $3 \leq k \leq n$.
\end{theorem}

\begin{theorem} {\rm \cite{KSW}}
There exists an $SMA(m, n; n, m)$ if and only if one of the following cases occurs:
\begin{itemize}
\item[$(1)$] $m = n = 1$;
\item[$(2)$] $m = 2$ and $n \equiv 0, 3 \pmod 4$;
\item[$(3)$]  $n = 2$ and $m \equiv 0, 3 \pmod 4$;
\item[$(4)$] $m,n > 2$.
\end{itemize}
\end{theorem}

Also the cases where each column contains 2 or 3 filled cells have been solved.
\begin{theorem}{\rm \cite{KE}}
There exists an $SMA(m, n; h, 2)$ if and only if one of the following cases occurs:
\begin{itemize}
\item[$(1)$] $m = 2$ and $n = h \equiv 0, 3 \pmod 4$;
\item[$(2)$] $m, h> 2$ and $mh = 2n$.
\end{itemize}
\end{theorem}

\begin{theorem}{\rm \cite{KLE}}
There exists an $SMA(m, n; h, 3)$ if and only if $3 \leq m$, $h\leq n$ and $mh = 3n$.
\end{theorem}

Necessary and sufficient conditions have been established also in the case when $h$ and $k$ are both even.

\begin{theorem}{\rm \cite{MP2}}
Let $h, k$ be two even integers with $h, k \geq 4$. There exists an $SMA(m, n; h, k)$ if and only if
$4\leq h  \leq n$, $4 \leq k \leq m$ and $mh = nk$.
\end{theorem}

Recent partial results for the general case can be found in \cite{MP3}. We summarize them in the following two theorems.

\begin{theorem}{\rm \cite{MP3}}
Let $m, n, h, k$ be four integers such that $3 \leq h \leq n$, $3 \leq k \leq m$ and $mh = nk$.
There exists an $SMA(m, n; h, k)$ whenever $\gcd(h, k) \geq 2$.
\end{theorem}

\begin{theorem}{\rm \cite{MP3}}
Let $m, n, h, k$ be four integers such that $3 \leq h \leq n$, $3 \leq k \leq m$ and $mh = nk$.
If $h\equiv 0 \pmod 4$, then there exists an $SMA(m, n; h, k)$.
\end{theorem}

\subsection{Non-zero sum Heffter arrays}

While studying the conjecture of Alspach (Conjecture \ref{Conj:als}), the authors of \cite{CDFP} proposed the following variant of a Heffter array which was termed a {\em non-zero sum Heffter array}.  The formal definition follows.

\begin{definition}\label{def:NZS}
Let $v=\frac{2nk}{\lambda}+t$ be a positive integer,
where $t$ divides $\frac{2nk}{\lambda}$,  and
let $J$ be the subgroup of $\Z_{v}$ of order $t$.
 A  $\lambda$-\emph{fold non-zero sum Heffter array $A$ over $\Z_{v}$ relative to $J$}, denoted by $^\lambda NH_t(m,n; h,k)$, is an $m\times n$ array
 with elements in $\Z_{v}$ such that:
\begin{itemize}
\item[$(\rm{a_4})$] each row contains $h$ filled cells and each column contains $k$ filled cells;
\item[$(\rm{b_4})$] the multiset $\{\pm x \mid x \in A\}$ contains each element of $\Z_v\setminus J$ exactly $\lambda$ times;
\item[$(\rm{c_4})$] the sum of the elements in every row and column {\underline {is different from $0$ }} in $\Z_v$.
\end{itemize}
\end{definition}
If the array is square, then $h=k\geq 1$, while if $m\neq n$ then at least one of $h$ and $k$ has to be greater than $1$.
A square non-zero sum Heffter array is denoted by $^\lambda NH_t(n;k)$. If $\lambda=1$ or $t=1$, then it is omitted.

\begin{example}\label{ex:NH75}
An $NH(7;5)$ over $\Z_{71}$:
\begin{center}
\begin{footnotesize}
$\begin{array}{|r|r|r|r|r|r|r|}\hline
10 &  & 16  & -1 & -2 & -3 &  \\
\hline   & 4 &  & -6 & -7 & -5 & 22   \\
\hline  -30 & 29 & 9 & -8 &   &  & 18  \\
\hline  -11 &   & -12 & -28 & -31 & 26 &  \\
\hline   & -14 & -15 & -13 &  & -17 & 25 \\
\hline  27 & -34 & 20 &   & -19 &  & -32  \\
\hline  24 & 23 &  &  & 21 &  -35 & 33  \\\hline
\end{array}$
\end{footnotesize}
\end{center}
\end{example}

In \cite{CDFP} the authors focused on the case in which $\lambda=1$ and they proved that the trivial necessary conditions are also sufficient in the very general case
in which the array is rectangular, empty cells are allowed and for every arbitrary subgroup $J$.   We would like to note, however, that in the {\em square case} the proof of existence is quite easy when there exists an $H_t(n; k)$.  Begin with the $H_t(n; k)$ and find a transversal of filled cells $T$ (this must exist). Then change the sign of the element in each of the cells of $T$ and note that now no row or column can add to $0$.
Hence we have constructed an $NH_t(n; k)$.  The general case is more involved and we emphasize that the proof
of the general result is not constructive.

\begin{theorem}{\rm \cite{CDFP}}\ 
There exists an $NH_t(m,n;h,k)$ if and only if $mh=nk$, $m\geq k \geq 1$, $n \geq h \geq 1$ and $t$ divides $2nk$.
\end{theorem}

These arrays have been introduced since, as with classical ones, they are related to open problems on partial sums, difference families, graph decompositions and biembeddings. The interested reader is referred to  \cite{CDFP} where  all the following results can be found.

By reasoning similar to what was done in Section \ref{sec:DF}, it is easy to see that if a row (column, respectively) of a $NH_t(m,n;h,k)$
admits an ordering whose partial sums are pairwise distinct and non-zero then one can construct a path of length $h$ ($k$, respectively).
For example, starting from the first row of the $NH(7;5)$ of Example  \ref{ex:NH75} and considering its
natural ordering one obtains the path $\G=[0,10,26,25,23,20]$
of length $5$, clearly $\Delta \G=\pm R_1$.
A non-zero sum Heffter array is said to be {\em simple} if every row and every column admits an ordering such that the partial sums are pairwise distinct and non-zero.
Hence, since every row and column of a non-zero sum Heffter array can be viewed as a set $T$ considered in Conjecture \ref{Conj:als}, we have that if
Conjecture \ref{Conj:als} were true, then every non-zero sum Heffter array would be simple. Investigating arrays having this additional property is important in view of the following results.

\begin{proposition}{\rm \cite{CDFP}}
If there exists a simple $^\lambda NH_t(m,n;h,k)$, then there exist
a $(2mh+t,t,P_h,\lambda)$-DF and a $(2nk+t,t, P_k,\lambda)$-DF.
\end{proposition}

\begin{proposition}{\rm \cite{CDFP}}
If there exists a simple  $^\lambda NH_t(m,n;h,k)$, then there exist
 a cyclic $P_h$-decomposition $\mathcal{D}$  and
 a cyclic $P_k$-decomposition $\mathcal{D}'$ both of $^\lambda K_{\frac{2nk+t}{t}\times t}$.
If $\lambda=1$ the decompositions $\mathcal{D}$ and $\mathcal{D}'$ are orthogonal.
\end{proposition}


In view of previous propositions, in \cite{CDFP,MePa} the authors focused on globally simple non-zero sum Heffter arrays,
obtaining the following existence results. All the proofs are constructive.

\begin{theorem}{\rm \cite{CDFP}}\label{NH1}
For every $n\geq k\geq 1$ there exists a  globally simple cyclically $k$-diagonal $NH(n;k)$.
\end{theorem}

\begin{theorem}{\rm \cite{CDFP}}
For every $m,n\geq 1$ there exists a globally simple $NH(m,n;n,m)$.
\end{theorem}

\begin{theorem}{\rm \cite{MePa}}\label{NH:t}
For every odd integer $n \geq 1$ and for every divisor $t$ of $n$, there
exists a globally simple $NH_t(n; n)$.
\end{theorem}

\begin{proposition}{\rm \cite{MePa}}\label{NH:t2}
For every odd integer $n \geq 1$ and for every $t\in\{2,2n,n^2,2n^2\}$, there
exists a globally simple $NH_t(n; n)$.
\end{proposition}

By Theorems \ref{NH1} and \ref{NH:t} and Proposition \ref{NH:t2} one gets the following result
(the case $n=2$ is trivial).

\begin{corollary}{\rm \cite{MePa}}
Let $n$ be a prime. There exists a globally simple $NH_t(n; n)$ for
every admissible $t$.
\end{corollary}

The case of $\lambda$-fold non-zero sum Heffter arrays with $\lambda>1$ has been investigated in \cite{CDFlambda}.
Here the authors consider a generic finite group $G$ and prove that
there exists a non-zero sum $\lambda$-fold Heffter array over $G$ relative to $J$
whenever the trivial necessary conditions are satisfied and $|G|\geq 41$.
Moreover, this value
can be decreased to $29$ when the array does not contain empty cells.

We  conclude this subsection with the theorem that gives the connection to biembeddings.

\begin{theorem}{\rm \cite{CDFP}}
Let $A$ be a non-zero sum Heffter array $NH(m,n;h,k)$ with  compatible orderings $\omega_r$
and $\omega_c$.
Then there exists a cellular biembedding of two circuit decompositions
 of $K_{2nk+1}$ into an orientable surface, such that the faces are multiples of $h$ and $k$ strictly larger than $h$
and $k$, respectively.
\end{theorem}

Details on  the faces of the biembeddings can be given considering the order of the row and column sums
in the group, as done in \cite{C,MePa}.

\subsection{Heffter configurations}

Heffter arrays are generalized to {\it Heffter configurations} in \cite{BP22}.
We recall that a $(v,b,k,r)$-configuration is a pair $(V,{\cal B})$ where $V$ is a set of $v$ {\it points} and
$\cal B$ is a $b$-set of $k$-subsets ({\it blocks}) of $V$ such that any two distinct points are contained
together in at most one block and any point occurs in exactly $r$ blocks.
Such a configuration is {\it resolvable} if there exists a partition of $\cal B$ into $r$ {\it parallel classes} each
of which is a partition of $V$, see \cite{BS21,Ge}.
Also, it is $\Z_v$-{\it additive}
if its point set is a subset of $\Z_v$ and each block sums to zero.

\begin{definition}
A $(nk,nr,k,r)$ \emph{Heffter configuration} (briefly HC) is a $\Z_{2nk+1}$-additive resolvable
$(nk,nr,k,r)$-configuration where the point set is a half-set of $\Z_{2nk+1}$.
It is {\it proper} if $r\geq3$ and  it is {\it simple} if each block admits a simple ordering.
\end{definition}

The above terminology is justified by the fact
that a  $(nk,nr,k,r)$-HC is equivalent to $r$ mutually sub-orthogonal Heffter systems  on the same half-set
of $\Z_{2nk+1}$.
Indeed each parallel class is nothing but a Heffter system D$(nk,k)$.
Hence,
a Heffter system D$(nk,k)$ is a $(nk,n,k,1)$-HC and an
 Heffter array $H(n;k)$, say $A$, is  a $(nk,2n,k,2)$-HC whose two parallel classes are the rows of $A$ and the columns of $A$.

The construction of proper HCs appears to be very hard. On the other
hand some infinite series have been obtained recursively by means
of small examples such as the following.

\begin{example}
The integer shiftable $H(5;4)$ below
gives a $(20,15,4,3)$-HC whose three parallel classes are:
 the rows;
 the columns;
and the diagonals $D_1, \ldots, D_5$.
$$
\begin{array}{|r|r|r|r|r|} \hline
-1& 2 & 17 & -18 & \\ \hline
6 & 10 & -13 & & -3\\ \hline
9 & -5 & & -8 & 4\\ \hline
-14 & & -16 & 11 & 19\\ \hline
 & -7 & 12 & 15 &-20 \\ \hline
\end{array}
$$
\end{example}

 Generalizing Proposition \ref{orth.systems}, we can state the following.
 \begin{theorem}
 A simple $(nk,nr,k,r)$-HC gives rise to $r$ mutually orthogonal cyclic $C_k$-decompositions of $K_{2nk+1}$.
 \end{theorem}

\subsection{Further generalizations}

Another natural generalization  is to work in an arbitrary additive group $G$ and not necessarily in the cyclic group.
Some work in this direction has been done in \cite{CPPBiembeddings} where the authors proposed the very
general concept of an \emph{Archdeacon array}.

\begin{definition}{\rm \cite{CPPBiembeddings}}
 An \emph{Archdeacon array} over an abelian group $(G, +)$ is an $m \times n$ array $A$
with elements in $G$, such that:
\begin{itemize}
  \item[$(\rm{b_5})$] the multiset $\{\pm x \mid x \in A\}$ contains each element of $G \setminus \{0\}$ at most once;
  \item[$(\rm{c_5})$] the elements in every row and column sum to $0$ in $G$.
\end{itemize}
\end{definition}

\begin{example}
An Archdeacon array on $\Z_{41}\oplus \Z_d$, for $d\geq 3$:
$$
\begin{array}{|r|r|r|r|r|} \hline
(-1,1) & (0,-1) & (-14,0) & (9,0) & (6,0) \\ \hline
(2,-1) & (-7,1) &  & (-5,0) & (10,0)\\ \hline
(17,0) & (12,0) & (-16,0) &  & (-13,0)\\ \hline
(-18,0) & (15,0) & (11,0) & (-8,0) & \\ \hline
 & (-20,0) & (19,0) & (4,0) & (-3,0) \\ \hline
\end{array}
$$
\end{example}

Clearly, Heffter arrays are a particular kind of Archdeacon arrays.
Examples, existence results and applications of these arrays can be found in Section $5$ of \cite{CPPBiembeddings}.

Another general concept termed a {\em quasi-Heffter array}  has recently been proposed in \cite{C}. In that paper the reader can find all the details
on the applications of these arrays to biembeddings and their full automorphism group.

\begin{definition}{\rm \cite{C}}
Let $v = 2nk + t$ be a positive integer, where $t$ divides $2nk$, and let $J$
be the subgroup of $\Z_v$ of order $t$. A \emph{quasi-Heffter array} $A$ \emph{over $\Z_v$ relative to $J$}, denoted
by $QH_t(m, n; h, k)$, is an $m \times n$ array with elements in $\Z_v$ such that:
\begin{itemize}
  \item[$(\rm{a_6})$] each row contains $h$ filled cells and each column contains $k$ filled cells;
 \item[$(\rm{b_6})$] the multiset $\{\pm x \mid x \in A\}$ contains each element of $\Z_v \setminus J$ exactly once.
 \end{itemize}
\end{definition}

So in this case the definition does not require any property on the row and column sums.
In general, given $S \subseteq \Z_v$ one may require that the sums of the rows and of the columns belong to $S$.
Clearly, if $S=\{0\}$ we find the classical concept of a Heffter array and  if $S=\Z_v\setminus\{0\}$ we get
the notion of a non-zero sum Heffter array.  For other choices of $S$ we have further variants of the classical concept having again applications to difference families, graph decompositions and biembeddings.

\section{Conclusions and Open Problems}\label{sec:conclusions}

Heffter arrays are a relatively new type of combinatorial design with original paper on this topic  published in 2015.  Since that time there have been a number of papers published on this subject as well as on topics that are related to Heffter arrays.  Here we have surveyed all the major results of these papers.
We believe that there are many more interesting questions relating to Heffter arrays and related subjects and we fully expect researchers to make good progress on these in the years to come.
We conclude this survey with a short list of open problems.
\begin{itemize}
\item Complete the construction of classical Heffter arrays in the rectangular case with empty cells, for known results see Theorems \ref{Fiore1} and \ref{Fiore2}.
\item Construct new infinite classes of globally simple Heffter arrays, for known results see Theorems \ref{thm:SH}, \ref{thm:globally_simple2} and \ref{thm:SH1}.
\item Determine new solutions to the Crazy Knight's Tour Problem, see Section \ref{sec:biemb}.
\item Determine new results concerning Conjectures \ref{Conj:als}, \ref{Conj:ADMS} and \ref{Conj:nostra}.
\item Establish more restrictive necessary conditions for integer $^\lambda H_t(n;k)$, in order to improve Proposition \ref{prop:necc}.
\end{itemize}

\subsection{Acknowledgements}
The authors would like to thank the Fields Institute for supporting the Stinson66 Conference and for financial support to the authors to attend the conference.

\end{document}